\documentclass[lettersize,journal]{IEEEtran}

\usepackage{amsmath,amsfonts, amssymb}
\usepackage{algorithmic}
\usepackage{algorithm}
\usepackage{array}
\usepackage[caption=false,font=normalsize,labelfont=sf,textfont=sf]{subfig}
\usepackage{textcomp}
\usepackage{stfloats}
\usepackage{url}
\usepackage{verbatim}
\usepackage{graphicx}
\usepackage{xcolor}
\usepackage{cite}
\usepackage{booktabs}
\usepackage{etoolbox}

\begin{document}


\newcommand{\varcolor}{black}
\newcommand{\varcolormoving}{black}

\newcommand{\n}{\mathrm{n}_{\config}}
\newcommand{\Zabcn}{\boldsymbol{Z}_{l}^{Car}}
\newcommand{\Rabcn}{\boldsymbol{R}_{l}^{Car}}
\newcommand{\Xabcn}{\boldsymbol{X}_{l}^{Car}}
\newcommand{\GMR}{\mathrm{GMR}_{\lineid}}
\newcommand{\rac}{{R}^{\text{ac}}_{\lineid}}

\newcommand{\Zabc}{\boldsymbol{Z}_{l}^{\text{Kr}}}
\newcommand{\Rabc}{\boldsymbol{R}_{l}^{\text{Kr}}}
\newcommand{\Xabc}{\boldsymbol{X}_{l}^{\text{Kr}}}

\newcommand{\Zij}{{\boldsymbol{Z}_{l,ij}}}
\newcommand{\Rij}{{\mathrm{\mathbf{R}}_{l,ij}}}
\newcommand{\Xij}{{\mathrm{\mathbf{X}}_{l,ij}}}

\newcommand{\Zin}{{\boldsymbol{Z}_{l,in}}}
\newcommand{\Rin}{{\mathrm{\mathbf{R}}_{l,in}}}
\newcommand{\Xin}{{\mathrm{\mathbf{X}}_{l,in}}}

\newcommand{\Znj}{{\boldsymbol{Z}_{l,nj}}}
\newcommand{\Rnj}{{\mathrm{\mathbf{R}}_{l,nj}}}
\newcommand{\Xnj}{{\mathrm{\mathbf{X}}_{l,nj}}}

\newcommand{\Znn}{{Z}_{l,nn}}
\newcommand{\Rnn}{{R}_{l,nn}}
\newcommand{\Xnn}{{X}_{l,nn}}

\newcommand{\A}{\boldsymbol{A}}
\newcommand{\Are}{\boldsymbol{A^{\text{re}}}}
\newcommand{\Aim}{\boldsymbol{A^{\text{im}}}}

\newcommand{\Ainv}{\boldsymbol{A}^{-1}}
\newcommand{\Ainvre}{\frac{1}{3}\Are}
\newcommand{\Ainvim}{\frac{-1}{3}\Aim}
\newcommand{\Ainvimnegative}{\frac{1}{3}\Aim}

\newcommand{\Zseq}{\boldsymbol{Z}_{l}^{012}}
\newcommand{\Rseq}{\boldsymbol{R}_{l}^{012}}
\newcommand{\Xseq}{\boldsymbol{X}_{l}^{012}}
\newcommand{\Yseq}{\boldsymbol{Y}^{\text{sh,012}}_{l}}

\newcommand{\Pabcn}{\boldsymbol{P}_{l}^{Car}}
\newcommand{\Pij}{{\boldsymbol{P}_{ij}}}
\newcommand{\Pin}{{\boldsymbol{P}_{in}}}
\newcommand{\Pnj}{{\boldsymbol{P}_{nj}}}
\newcommand{\Pnn}{\symbol{P}_{l}^{nn}}

\newcommand{\Cabcn}{\boldsymbol{C}^{\text{sh}}_{l}}
\newcommand{\Yabcn}{\boldsymbol{Y}^{\text{sh}}_{l}}
\newcommand{\Yij}{{\boldsymbol{Y}^{\text{sh}}_{l,ij}}}
\newcommand{\Yin}{{\boldsymbol{Y}^{\text{sh}}_{l,in}}}
\newcommand{\Ynj}{{\boldsymbol{Y}^{\text{sh}}_{l,nj}}}
\newcommand{\Ynn}{Y^{\text{sh}}_{l,nn}}
\newcommand{\Yabc}{\boldsymbol{Y}^{\text{sh,abc}}_{l}}

\newcommand{\Sij}{{S}_{\config,ij}}
\newcommand{\Sii}{{S}_{\config,ii}}
\newcommand{\R}{{R}_{\lineid}}

\newcommand{\Kgmr}{{K}^{\text{gmr}}_{\config}}

\newcommand{\radius}{{r}_{\lineid}}
\newcommand{\x}{\boldsymbol{x_{\config}}}
\newcommand{\y}{\boldsymbol{y_{\config}}}
\newcommand{\Dij}{{D}_{\config,ij}}
\newcommand{\Dii}{{D}_{\config,ii}}
\newcommand{\Kradius}{{K}^{\text{r}}_{\config}}
\newcommand{\tnom}{{t}^{\text{nom}}}
\newcommand{\Rnom}{{R}^{\text{nom}}_{\lineid}}

\newcommand{\freevar}[2]{{#1}_{#2}}
\newcommand{\vref}{v_{\config}^{\text{ref}}}
\newcommand{\vrefmin}{v_{\config}^{\text{ref,min}}}
\newcommand{\vrefmax}{v_{\config}^{\text{ref,max}}}

\newcommand{\Dmin}{D_{\config}^{\text{min}}}
\newcommand{\uminOH}{u_{\config}^{\text{min,OH}}}
\newcommand{\umaxOH}{u_{\config}^{\text{max,OH}}}
\newcommand{\umincable}{u_{\config}^{\text{min,cable}}}
\newcommand{\umaxcable}{u_{\config}^{\text{max,cable}}}

\newcommand{\Rzero}{R_{\lineid,00}}
\newcommand{\Rone}{R_{\lineid,11}}
\newcommand{\Xzero}{X_{\lineid,00}}
\newcommand{\Xone}{X_{\lineid,11}}

\newcommand{\Rzeroref}{\Rzero^{\text{ref}}}
\newcommand{\Roneref}{\Rone^{\text{ref}}}
\newcommand{\Xzeroref}{\Xzero^{\text{ref}}}
\newcommand{\Xoneref}{\Xone^{\text{ref}}}

\newcommand{\Rzeroaux}{\Rzero^{\text{aux}}}
\newcommand{\Roneaux}{\Rone^{\text{aux}}}
\newcommand{\Xzeroaux}{\Xzero^{\text{aux}}}
\newcommand{\Xoneaux}{\Xone^{\text{aux}}}

\newcommand{\Rzerostar}{\Rzero^{*}}
\newcommand{\Ronestar}{\Rone^{*}}
\newcommand{\Xzerostar}{\Xzero^{*}}
\newcommand{\Xonestar}{\Xone^{*}}

\newcommand{\Bzero}{B_{\lineid,00}}
\newcommand{\Bone}{B_{\lineid,11}}

\newcommand{\Bzeroref}{\Bzero^{\text{ref}}}
\newcommand{\Boneref}{\Bone^{\text{ref}}}

\newcommand{\Bzeroaux}{\Bzero^{\text{aux}}}
\newcommand{\Boneaux}{\Bone^{\text{aux}}}

\newcommand{\Bzerostar}{\Bzero^{*}}
\newcommand{\Bonestar}{\Bone^{*}}

\newcommand{\p}{\rho_{\material}}
\newcommand{\T}{T_{\lineid}}
\newcommand{\Tcoef}{\alpha_{\material}}
\newcommand{\Cskin}{C^{s}_{\lineid}}
\newcommand{\Cprox}{C^{p}_{\lineid}}
\newcommand{\rdc}{{R}^{dc}_{\lineid}}
\newcommand{\N}{{N}_{\config}}
\newcommand{\area}{{A}_{\lineid}}

\newcommand{\rdcACSR}{\boldsymbol{R^{dc,ACSR}_i}}
\newcommand{\rdcAl}{\boldsymbol{R^{dc,Al}_i}}
\newcommand{\rdcSt}{\boldsymbol{R^{dc,St}_i}}
\newcommand{\racACSR}{\mathrm{R}^{ac,ACSR}_{i}}
\newcommand{\CACSR}{\mathrm{C}^{ACSR}_{i}}

\newcommand{\Zdiff}{{Z}_l^{\text{diff,series}}}
\newcommand{\BZdiff}{{Z}_l^{\text{diff}}} 
\newcommand{\Zthre}{{Z}^{\text{thre}}}
\newcommand{\Zref}[1]{{Z_{l,#1}}^{\text{ref}}}

\newcommand{\config}{f}
\newcommand{\lineid}{l}
\newcommand{\material}{m}
\newcommand{\lineconfig}{q}
\newcommand{\combination}{\lineid \config \material}
\newcommand{\var}{\varphi}

\newcommand{\setconfig}{\mathcal{F}}
\newcommand{\setlineid}{\mathcal{L}}
\newcommand{\setmaterial}{\mathcal{M}}
\newcommand{\setlineconfig}{\mathcal{Q}}
\newcommand{\setcombination}{\mathcal{G}}
\newcommand{\setvar}{\mathcal{\phi}}

\newcommand{\Tmin}{\T^{\text{min}}}
\newcommand{\Tmax}{\T^{\text{max}}}
\newcommand{\rmin}{\radius^{\text{min}}}
\newcommand{\rmax}{\radius^{\text{max}}}
\newcommand{\Amin}{\area^{\text{min}}}
\newcommand{\Amax}{\area^{\text{max}}}

\newcommand{\slack}{\beta}
\newcommand{\frequency}{f^{\text{fund}}_l }

\newcommand{\imag}{\textcolor{\varcolor}{\texttt{j}}}
\newcommand{\tick}{\checkmark}
\newcommand{\OHconfig}{\config^\text{OH}}
\newcommand{\Cableconfig}{\config^\text{cable}}
\newcommand{\OHstdconfig}{\config^{l,\text{OH}}}
\newcommand{\setOHstdconfig}{\setconfig^{l,\text{OH}}}
\newcommand{\setOHconfig}{\setconfig^{l,\text{OH}}}
\newcommand{\setcableconfig}{\setconfig^{l,\text{cable}}}

\newcommand{\AminOHfor}{\A^\text{min,OH}}
\newcommand{\AmaxOHfor}{\A^\text{min,OH}}
\newcommand{\TminOHfor}{\T^\text{min,OH}}
\newcommand{\TmaxOHfor}{\T^\text{min,OH}}

\newcommand{\AminOHCable}{\A^\text{min,cable}}
\newcommand{\AmaxOHCable}{\A^\text{min,cable}}
\newcommand{\TminOHCable}{\T^\text{min,cable}}
\newcommand{\TmaxOHCable}{\T^\text{min,cable}}

\newcommand{\ua}{\freevar{u}{\config,1}}
\newcommand{\ub}{\freevar{u}{\config,2}}
\newcommand{\va}{\freevar{v}{\config,1}}

\newcommand{\uastd}{\ua^\text{std}}
\newcommand{\ubstd}{\ub^\text{std}}
\newcommand{\vastd}{\va^\text{std}}
\newcommand{\vrefstd}{v_{\config}^{\text{ref,std}}}
\newcommand{\rstd}{\radius^\text{std}}
\newcommand{\Astd}{\area^\text{std}}

\newcommand{\kone}{k_1^{\text{Car}}}
\newcommand{\ktwo}{k_2^{\text{Car}}}
\newcommand{\kthree}{k_3^{\text{Car}}}
\newcommand{\kfour}{k_4^{\text{Car}}}
\newcommand{\kfive}{k_5^{\text{Car}}}

\newcommand{\tnommin}{{t}^{\text{nom,min}}}
\newcommand{\tnommax}{{t}^{\text{nom,max}}}

\title{The Inverse Carson's Equations Problem: Definition, Implementation and \textcolor{\varcolor}{Numerical} Experiments}
\author{C. H. Tam, F. Geth, N. Mithulananthan, ~\IEEEmembership{Senior Member,~IEEE,}}
\markboth{Submitted to IEEE TRANSACTIONS ON POWER DELIVERY, VOL. xx, NO. x, April 2024}%
{Shell \MakeLowercase{\textit{et al.}}: A Sample Article Using IEEEtran.cls for IEEE Journals}


\maketitle

\begin{abstract}
In recent years, with the increase in renewable energy and storage penetration, power flow studies in low-voltage networks have become of interest in both industry and academia. 
Many studies use impedance represented by sequence components due to the lack of datasets with fully parameterized impedance matrices. 
This assumes that the network impedance is balanced, which is typically not the case in the low-voltage network and therefore risks the accuracy of the study. 
This paper proposes a methodology for the recovery of more detailed impedance data from sequence components as an inverse problem, i.e. the inverse Carson's equations problem, for both overhead lines and cables. 
We consider discrete properties like material and configuration of conductors common in the distribution network and investigate what data can be reliably recovered from only sequence components using nonlinear optimisation models. Presented results include uniqueness of recovered variables and the likelihood of mismatch. 
\end{abstract}
\textcolor{\varcolor}{
\section*{NOMENCLATURE}\label{Nomenclature}}
\noindent
\textcolor{\varcolor}{
\textit{Conductor Variables} \vspace{3pt}\\
\begin{tabular}{@{} p{2.1cm} p{6cm}} 
$\radius, \area $ & Conductor strand radius and area.\\
$\T $ & Conductor temperature.\\
$\rdc,\rac $ & Conductor dc and ac resistance. \\ 
$\GMR$ & Geometric mean radius of a conductor. \\
$\tnom $ & Conductor insulation thickness. \\
$\R,\Rnom$ & Conductor radius, with insulation.  \\
$\ua,\ub $ & Geometry variables for x-coordinate. \\
$\va, \vref $ & Geometry variables for y-coordinate.  \\
$\x,\y$ & Conductor x- and y-coordinate vectors.\\
$\Dij,\Sij $ & Distance and image distance matrix. \\
\end{tabular}}
\textcolor{\varcolor}{
\textit{Impedance \& admittance variables} \vspace{3pt}\\
\begin{tabular}{@{} p{2.1cm} p{6cm}}
$\Zabcn $ & Carson's derived impedance matrix. \\
$\Zabc $ & Kron's reduced impedance matrix.\\
$\Zseq $ & Sequence impedance matrix.\\
$\Pabcn,\Cabcn  $ & Potential coefficient, capacitance matrix. \\
$\Yabcn,\Yseq $ & Shunt admittance, sequence matrix.\\
$\Zdiff $ & Error over given diag. seq. impedance. \\
$\BZdiff $ & Above with given diag. seq. admittance. \\
\end{tabular}}
\textcolor{\varcolor}{
\textit{Bound} \vspace{3pt}\\
\begin{tabular}{@{} p{2.1cm} p{6cm}} 
$\rmin,\rmax $ & Min. \& max. conductor strand radius.\\
$\Amin,\Amax $ & Min. \& max. conductor area.\\
$\Tmin,\Tmax $ & Min. \& max. conductor temperature.\\
$\tnommin,\tnommax $ & Min. and max. of $\tnom$. \\
$\vrefmin,\vrefmax $ & Min. \& max. conductor reference height. \\
$\Dmin$ & Min. distance between overhead wires. \\
$\umaxOH $ & Half-side crossarm length of pole. \\
$\umincable $ & Min. distance between cable cores. \\
$\umaxcable $ & Max. distance between cable cores. \\
\end{tabular}}
\textcolor{\varcolor}{
\textit{Parameter} \vspace{3pt} \\
\begin{tabular}{@{} p{2.1cm} p{6cm}} 
$\imag$ & Imaginary unit. \\
$\frequency $ & Fundamental frequency of power grid. \\
$k^{\text{Car}} $ & Carson's equations constants. \\
$\A$ & Symmetrical component transf. matrix.\\
$\n $ & Number of conductors of a line. \\
$\N$ & Number of strands of a conductor. \\
$\Kgmr, \Kradius$ & Ratio between $\radius$ \& $\GMR$, $\radius$ \& $\R$. \\
$\p,\Tcoef $ & Material resistivity, temperature coef. \\
$\Cskin,\Cprox $ & Skin \& proximity effect correction factor.\\
$\rstd $ & Standardised conductor strand radius. \\
$\uastd,\ubstd,\vastd$ & Standardised geometry parameters. \\
$\vrefstd$ & Standardised reference height. \\
$\Bzeroref, \Boneref $ & Given zero \& pos. sequence susceptance.\\
$\Rzeroref, \Roneref$ & Given zero \& pos. sequence resistance.\\
$\Xzeroref, \Xoneref $ & Given zero \& pos. sequence reactance.\\
$\slack $ & Slack for given sequence components.\\
\end{tabular}}
\textcolor{\varcolor}{
\textit{Set and indices} \vspace{3pt} \\
\begin{tabular}{@{} p{2.1cm} p{6cm}}
Line  &    $\lineid \in \setlineid $\\
Material   &  $\material \in \setmaterial =\{\text{Al-1350},\text{Cu}\}$\\
Configuration & $\config \in \setconfig$ \\
Combination & $ \lineid \config \material \in \setcombination \subseteq \setlineid \times \setconfig \times  \setmaterial$\\  
Conductor  &    $i,j \in \mathcal{W} = \{a,b,c,n \} $\\
Phase &    $p \in \mathcal{P} = \{a,b,c \}$ \\
Decision var. & $\var \in \setvar$ \vspace{5pt}\\
\end{tabular}}
\\

\section{Introduction}\label{section_intro}
\IEEEPARstart{W}ITH the proliferation of roof-top photovoltaics, electric vehicles and battery storage systems, concerns have been raised on the operation of  (unbalanced) low-voltage (LV) distribution networks. Power flow simulation is a technique used by utilities to understand and manage the impact of these new energy technologies. One obstacle for LV power flow studies is the lack of high-quality network models, which nevertheless are needed for decision-support approaches\cite{beckstedde_fit_2023}.

\subsection{Imperfect Network Datasets}
\textcolor{\varcolor}{
Because LV networks have untransposed lines, the positive-sequence approximation for impedance does not hold. Historically, utilities may have calculated the impedance of untransposed lines (with or without neutral) for overhead (OH) lines and cables, using information readily available from datasheets such as ac resistance of conductors at some temperature (e.g. 75\,°C). Information for calculating impedance anew, such as conductor and geometry types may not be kept properly, making recovery of full impedance matrices difficult. Projects working on releasing LV network models, such as \cite{CSIRO_LV_report} in Australia and \cite{Manchester_LV_report} in the UK, do not include the more detailed wire-coordinate impedances but instead provide $Z_{11}$ and $Z_{00}$\footnote{\textcolor{\varcolor}{we refer to these as `diagonal sequence components'}}. This cannot represent the untransposed lines in LV networks, and can impact the accuracy of power flow simulation \cite{tam_inclusive_2022,kersting_whys_2011,claeys_applications_2021}.  
}
\textcolor{\varcolor}{
\subsection{Importance of Having Good Impedance Data  }
Recovering conductor and geometry information from diagonal sequence components has several advantages. Firstly, one can recompute the neutral-explicit impedance for better study of network congestion under normal operation. Unlike transmission networks, LV networks have unbalanced loading resulting in non-negligible neutral voltage. A maximum neutral voltage of 10\,V is part of the voltage supply standard in Australia \cite{Essential_Energy_supply_2014} and should be considered in studies. 
Distribution utilities examine network congestion through state estimation \cite{Energex_DSSE}, which relies on high-quality impedance data \cite{geth_data_2023}. Accuracy is important for decision-making as lots of investment is being made into grid development and upgrades\cite{Electricity_cost_AER}.
Lastly, recovering conductor type helps to derive other properties, such as its current limit for fault analysis.
}

\subsection{Carson's Equations}
Impedance of overhead lines and cables are dependent on the environment and installation, due to mutual inductive effects between conductors and capacitive coupling to earth. 
Impedance values are established from Maxwell's laws, either solved through finite element simulation or Carson's equations. 
\textcolor{\varcolor}{
Carson’s equations~\cite{carson_wave_1926}, derived by approximating
the Maxwell’s laws, obtain the impedance matrices for electromagnetically coupled conductors based on geometry, material, resistivity of the conductor and the earth, temperature, and cross-section, taking into account ground return current.
}
\subsection{Contribution}
\textcolor{\varcolor}{
To enhance the quality of impedance datasets, this paper proposes a methodology to recover the inputs to Carson's equations from zero and positive sequence components derived from untransposed lines. The contributions are:
\begin{itemize}
    \item Proposing a methodology to solve the inverse Carson's equations problem, to (partially) recover conductor and geometry information for overhead lines and cables.
    \item Recomputing the full unbalanced impedance of overhead lines and cable through recovered information. 
    \item Conducting numerical simulations to validate accuracy, model uncertainty and mismatch risks.
    \item Demonstrating recovery with missing data and showing the need for data validation, through utility case study.
\end{itemize}
}

\subsection{Notation and Preliminaries}
\begin{figure}[tbh]
\centering
\includegraphics[width=1\linewidth]{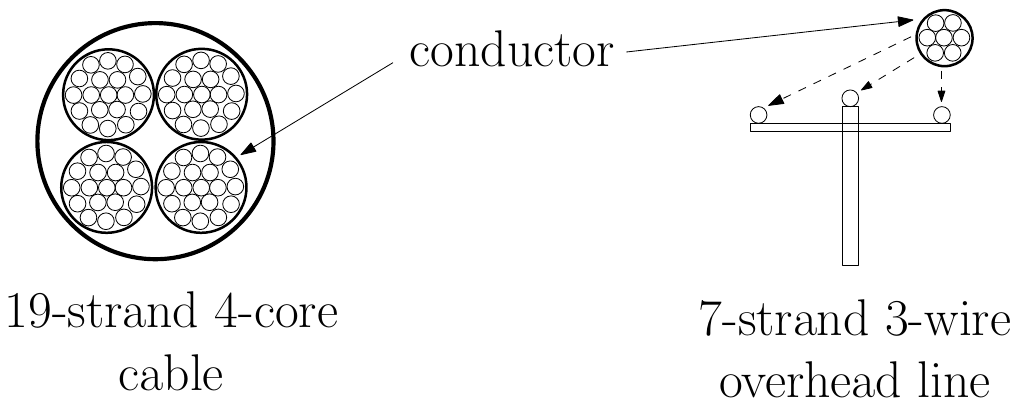}
\caption{The left diagram represents a cable with $\N$=19 and $\n$=4, on the right  an overhead line with $\N$=7 and $\n$=3.}
\label{fig:preliminary_figure}
\end{figure}

To prevent ambiguity, specific terms are defined to distinguish overhead line and cable construction properties. 
A \emph{conductor} is formed by $\N$ \emph{strands} and it is more specifically referred to as a \emph{wire} for \emph{overhead lines} and as a \emph{core} for \emph{cables}.
For example, on the left of Fig. \ref{fig:preliminary_figure} is a 4-core 19-strand cable and on the right is a 3-wire 7-strand overhead line.

The term \emph{forward calculation} is referred to as having all required information to obtain the reference diagonal sequence components of impedance $\Rzeroref$, $\Roneref$,$\Xzeroref$ $\Xoneref$, and that of shunt susceptance $\Bzeroref$, $\Boneref$. Next, the term \emph{inverse estimation} is referred to as establishing a likely combination to explain the sequence components being given. 

\begin{figure*}[th!]
\centering
\includegraphics[trim={0 0.1cm 0 0cm},clip,width=1\linewidth]{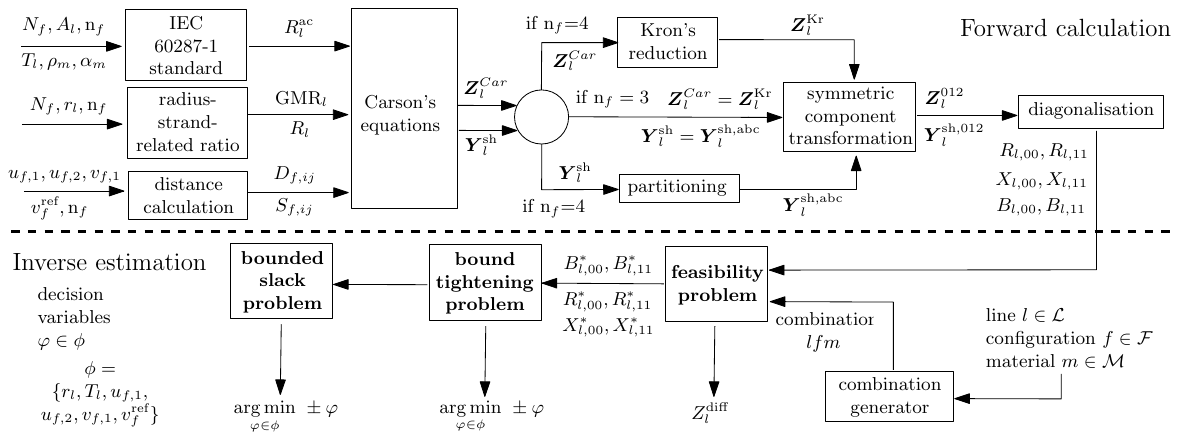}
\caption{Computational steps in the forward and inverse Carson's equations solution process.}
\label{fig:flow_block_diagram}
\end{figure*}

\textcolor{\varcolormoving}{In forward calculation as depicted in the top half of Fig. \ref{fig:flow_block_diagram}, the inputs to Carson's equations are taken from standard construction codes. For a 3- or 4-conductor line, sequence components are obtained. Carson's equations are used to obtain the series impedance and shunt admittance matrices. For a 4-conductor line, Kron's reduction of impedance matrix and partitioning of admittance matrix are done to reduce the matrix so that symmetrical component transformation can be performed to obtain the sequence components value.} 

\subsection{Paper Structure}

Section \ref{section: lit review} is a literature review on the importance of using accurate LV network impedance, and the use of Carson's equations in power systems. 
Section \ref{section: network modelling} summarizes how to derive the impedance values for common LV networks configurations through Carson's equations. 
Section \ref{sec_opt_prob_spec} defines the mathematical problem statement of the inverse Carson's equations problem, and 
Section \ref{section: test_case_and_validation} demonstrates recovery of line parameter through untransposed lines.
Lastly, Section \ref{sec_conclusions} provides conclusions and avenues for future work.

\section{Literature Study} 
\label{section: lit review}

\subsection{Consequence of Approximation in LV Network Impedance}
The perfect-grounding assumption is equivalent to Kron's reduction of the neutral in a 4-conductor network, whereas the diagonalisation step assumes the original 4-conductor network to be balanced as shown in Fig. \ref{fig:flow_block_diagram}. Urquhart \cite{urquhart_accuracy_2016} shows that assuming zero off-diagonal sequence components causes 17\% voltage error  for a 4-core cable of 100\,m in the UK even when a balanced current of 50\,A is used. 
Claeys et al. \cite{claeys_optimal_2022} compare  4-wire networks with explicit grounding and the equivalent Kron-reduced networks for distributed generator dispatch optimisation. They show that using the optimisation result from a Kron-reduced network can lead up to 4\% phase-to-neutral error when validated against a 4-wire model.  
Our work \cite{tam_inclusive_2022} found that significant underestimation of voltage violation occurs under unbalanced loading if network impedance is assumed balanced. Lastly, using correct impedance data is especially important in contexts such as distribution state estimation, as any model error (e.g. impedance error) will get mixed into the measurement residuals, thereby impacting the most likely network state estimate \cite{vanin_combined_2023}.

\subsection{Observed Accuracy of Carson's Equations}
Carson's equations model the earth as a homogeneous semi-definite solid, creating conductor image underneath the ground and deriving impedance fundamentally by waveform propagation of electrical signal \cite{kersting_series_2012,carson_wave_1926}. 
The original equations are further simplified to the \emph{modified} Carson's equations by eliminating image-dependent coordinate variables. 
The approximation made by the modified equations has an error of less than 1$\%$ compared with its original form \cite{kersting_series_2012}. 
Although Carson's equations are at first derived for overhead lines, they also work well for underground cables with error within 1\% compared with finite element analysis at fundamental frequency \cite{urquhart_series_2015}. 
It is shown to be insensitive to ground resistivity and cable depth values \cite{urquhart_series_2015}.
\textcolor{\varcolor}{Other than Carson’s equations, a variety of approximations also exists to simplify the infinite integral term associated with deriving the self and mutual impedance.} Ref. \cite{keshtkar_improving_2014}  details the derivation of impedance and compares  approximations of the infinite integral terms between the models from Deri, Alvarado, Noda, Pizarro, Dubanton and Carson, and shows that other methods mainly deviate from Carson's for long distribution lines. 
Supported by this evidence, the authors believe that the modified Carson's equations provide an attractive trade-off between accuracy and simplicity (i.e. closed form, algebraic) and is therefore selected as the method to establish impedances in LV networks.
Another approach to derive impedance is through finite element analysis (FEA),
where detailed schematics of cross-sections of cables, the material types and insulation medium involved are required. 
Boundary conditions for the solution space are needed for truncating the magnetic field \cite{urquhart_series_2015}. 
Because conductors and medium are divided into elements (the mesh), a nonuniform structure can be defined and evaluated. 
Because of its relative complexity, FEA is often used to validate the assumption of models such as  Carson's equations, rather than modelling all cables and lines in the network\cite{cristofolini_comparison_2017}. 

\subsection{Inverse Problems in Power Systems}
\textcolor{\varcolor}{
Inverse problems in network models focus on recovery of nodal admittance \cite{Yuan2023,low2024reverse} or network topology \cite{yizheng_liao_distribution_2015, deka_estimating_2016} from measurement data.
Recovery of untransposed impedances from diagonal sequence components is left unexplored. In LV networks, diagonal sequence components and line length are commonly provided separately, and may not be correctly given or validated. State-of-the-art includes work by Vanin et al. \cite{vanin_combined_2023}, who propose line length recovery considering unbalanced impedance in LV network, based on voltage magnitude and power measurement from smart meters. However, cable types are assumed known. Yuan et. al\cite{Yuan2023} uses graph theory to show that single-phase network admittance can be uniquely recovered with phasor measurement units, and when nodes without measurement are connected with at least 3 nodes with measurement. Low in \cite{low2024reverse} extends the work in \cite{Yuan2023} to multi-phase network and shows that admittance can be recovered when the network is radial and has the same line type. Instead of recovering impedance or admittance through measurements, we contribute by recovering properties of overhead line and cables from given sequence impedances. 
}
\section{Impedance Models for Low-Voltage Networks} \label{section: network modelling}
\subsection{Derivation of Impedance from First Principles}
We use $i,j \in \mathcal{W} = \{a,b,c,n \} $ to refer to different conductors of an overhead line or cable, the conductor can be a neutral conductor $n$ or a phase conductor with phase $p \in \mathcal{P} = \{a,b,c \}$. 
An overhead line or cable has three indices: the line id $\lineid \in \setlineid $, the material  $\material \in \setmaterial =\{\text{Al-1350},\text{Cu}\}$, and the spatial configuration $\config \in \setconfig$.


The per-length impedance matrix symbols for a 4-conductor overhead line and cable $l$ with phase and neutral conductors $abcn$ are defined, 
\begin{IEEEeqnarray}{c}
\Zabcn \!=\! \Rabcn \!+\! \imag \Xabcn\!=\! 
\left( \begin{array}{c c c | c} 
Z_{l,aa} \!&\! Z_{l,ab} \!&\! Z_{l,ac} \!&\! Z_{l,an} \\
Z_{l,ba} \!&\! Z_{l,bb} \!&\! Z_{l,bc} \!&\! Z_{l,bn} \\
Z_{l,ca} \!&\! Z_{l,cb} \!&\! Z_{l,cc} \!&\! Z_{l,cn} \\
\hline
Z_{l,na} \!&\! Z_{l,nb} \!&\! Z_{l,nc} \!&\! Z_{l,nn} \label{eq_Z_abcn} 
\end{array} \right)  \nonumber \\ 
\!=\! \left(
\begin{array}{c | c}
\Zij & \Zin \\
\hline 
\Znj &  \Znn 
\end{array}
\right) [\Omega /km].
\end{IEEEeqnarray}
\textcolor{\varcolor}{
The Carson's equations for per-length self and mutual impedance in SI unit are provided by Cleenwerck et al. in \cite{cleenwerck_approach_2022}.
The self impedance for wire $i \in \{a,b,c,n\}$ is,
\begin{IEEEeqnarray}{l}
Z_{l,ii}^{\text{Car}} = \rac \!+\! \kone \!+\! \imag \cdot \ktwo \left(ln\frac{1}{\kthree \cdot \GMR} + \kfour \right)\!, \label{Carson_self}
\end{IEEEeqnarray} 
and the mutual impedance between wires $i$ and $j$ is, 
\begin{IEEEeqnarray}{l}
Z_{l,ij}^{\text{Car}} = \kone  \!+\! \imag \cdot \ktwo \left( ln \frac{1}{\kthree \cdot \Dij} + \kfour \right) \label{Carson_mutual},
\end{IEEEeqnarray}
where 
$\rac$ [$\Omega/km$] is the ac resistance of a conductor,
$\GMR$ [mm] is the geometric mean radius of a conductor,
$\Dij$ [mm] is the distance between conductors. 
The remaining constants in Carson's equations are tabulated  in Table \ref{tab_Carson_constants}.
From \eqref{Carson_self}-\eqref{Carson_mutual}, it is observed that $\Zabcn$ is  symmetric in both the real and imaginary parts.
}
\begin{table}[tbh]
\textcolor{\varcolor}{
\caption{Carson's equation constants \cite{cleenwerck_approach_2022}.}
\label{tab_Carson_constants}
    \centering 
    \begin{tabular}{l l l l l l }
        \toprule
    Cons. & value & unit & Cons. & value & unit \\
        \midrule
    $\kone$ & 0.049348 & [$\Omega/km$] & $\ktwo$ &0.062832 & [$\Omega/km$]\\
    $\kthree$ & 3.28084$\cdot10^{-3}$ & [mm$^{-1}$] & $\kfour$ & 8.0252 & [--] \\
    $\kfive$ & 17.98742 & [km/uF] & & & \\
        \bottomrule
    \end{tabular}
    }
\end{table}

The ac resistance is a function $h$ of \textcolor{\varcolor}{conductor} area $\area$, temperature $\T$ and material $\material$. The material determines resistivity $\p$ and temperature-coefficient $\Tcoef$,
\begin{IEEEeqnarray}{l}
\rac= h(\p, \Tcoef, \T, \area ) \label{eq_rac_dummy}.
\end{IEEEeqnarray} 

Kron's reduction of the neutral is performed to transform the network into a reduced $3\times 3$ matrix $\Zabc$ by assuming neutral voltage to be zero, 
\begin{IEEEeqnarray}{c}
\Zabc= \Rabc+\imag\Xabc  =  \Zij - \frac{1}{\Znn}   \Zin  \Znj. \label{Kron reduction eq}
\end{IEEEeqnarray}
The sequence impedance matrix $\Zseq$ is now derived from the Kron-reduced matrix, 
\begin{IEEEeqnarray}{c}
\Zseq =\Rseq+\imag\Xseq =\Ainv \Zabc \A \label{SCT eq}.
\end{IEEEeqnarray} 
The entries in the sequence matrix are,
\begin{IEEEeqnarray}{c}
\Zseq =
\begin{pmatrix}
    Z_{l,00} & Z_{l,01} & Z_{l,02} \\
    Z_{l,10} & Z_{l,11} & Z_{l,12} \\
    Z_{l,20} & Z_{l,21} & Z_{l,22} 
\end{pmatrix}.\label{eq_Z012_definition}
\end{IEEEeqnarray} 
Note that the diagonal entries are called to zero (00), positive (11) and negative sequence (22) impedance, and the positive and negative impedance are equal.
Although $\Zabcn$ is symmetric, the sequence transformation generally destroys that symmetry, e.g. $Z_{l,10} \neq Z_{l,01}$.

The symmetrical component transformation matrix $\A$ is,
\begin{IEEEeqnarray}{c}
\A 
=  \Are + \imag \Aim 
= 
\begin{pmatrix}
    1 & 1 & 1 \\
    1 & \alpha^2 & \alpha \\
    1 & \alpha & \alpha^2 \\
\end{pmatrix} \nonumber \\
= 
\begin{pmatrix}
    1 & \phantom{-}1 & \phantom{-}1 \\
    1 & -\frac{1}{2} & -\frac{1}{2} \\
    1 & -\frac{1}{2} & -\frac{1}{2} \\
    \end{pmatrix}
+ \imag 
\begin{pmatrix}
    0 & \phantom{-}0 & \phantom{-}0 \\
    0 & -\frac{\sqrt3}{2} & \phantom{-}\frac{\sqrt3}{2} \\
    0 & \phantom{-}\frac{\sqrt3}{2} & -\frac{\sqrt3}{2} \\
    \end{pmatrix} 
    \nonumber ,
\end{IEEEeqnarray}
with $\alpha=e^{2\pi/3}$, and its inverse is,
\begin{IEEEeqnarray}{c}
\Ainv =  \frac{1}{3}\Are+ \imag\frac{-1}{3}\Aim \nonumber.
\end{IEEEeqnarray}

Despite shunt admittance being  small to be negligible for overhead lines, underground cables have a  higher shunt admittance\cite{kersting_shunt_2012}. 
Although unfortunately sequence shunt admittance is usually not recorded by utilities, the modelling of sequence shunt admittance is provided in this paper for completeness of modelling.
Shunt conductance due to leakage current between conductors and insulation medium is neglected because it is very small compared to shunt susceptance \cite{kersting_shunt_2012}. 

To derive sequence admittance for shunt, the potential coefficient matrix $\Pabcn$ of line $l$ is firstly derived \cite{kersting_shunt_2012},
\begin{align} 
 \Pabcn &=
\left( \begin{array}{c c c c}
P_{l,aa} & P_{l,ab}& P_{l,ac} & P_{l,an} \\
P_{l,ba} & P_{l,bb}& P_{l,bc} & P_{l,bn} \\
P_{l,ca} & P_{l,cb}& P_{l,cc} & P_{l,cn} \\
P_{l,na} & P_{l,nb}& P_{l,nc} & P_{l,nn} \label{eq_Pabcn} 
\end{array} \right)  \textcolor{\varcolor}{[km/uF]}.
\end{align}
\textcolor{\varcolor}{
The self and mutual entries are defined, $i,j \in \{a,b,c,n \}$,
\begin{IEEEeqnarray}{c}
P_{l,ii} = \kfive \cdot ln \left( \frac{\Sii}{\R}\right) \label{P_self},\IEEEyesnumber \IEEEyessubnumber\\
P_{l,ij} = \kfive \cdot ln \left( \frac{\Sij}{\Dij} \right),  \IEEEyessubnumber \label{P_mutual}
\end{IEEEeqnarray}
where $\Sii$ is the distance [$mm$] between conductor $i$ of line $l$ and its image below ground. Next, $\Sij$ [$mm$] is the distance between actual conductor $i$ and image of another conductor $j$. Variable $\R$ is the conductor overall radius [$mm$] without insulation. Lastly, $\kfive$ is a constant in Table \ref{tab_Carson_constants}.}

The capacitance matrix $\Cabcn$ is the matrix inverse of $\Pabcn$,
\begin{IEEEeqnarray}{c}
\Cabcn = (\Pabcn)^{-1} \label{C_equation},
\end{IEEEeqnarray} 
and is used to compute the shunt admittance matrix $\Yabcn$ at  fundamental frequency $\frequency$ , 
\begin{IEEEeqnarray}{c}
\Yabcn = j  2 \pi \frequency \Cabcn,
\end{IEEEeqnarray} 
which also gets partitioned,
\begin{IEEEeqnarray}{c}
\Yabcn=\left( \begin{array}{c c c | c}
Y_{l,aa} & Y_{l,ab}& Y_{l,ac} & Y_{l,an} \\
Y_{l,ba} & Y_{l,bb}& Y_{l,bc} & Y_{l,bn} \\
Y_{l,ca} & Y_{l,cb}& Y_{l,cc} & Y_{l,cn} \\
\hline
Y_{l,na} & Y_{l,nb}& Y_{l,nc} & Y_{l,nn} \label{first_equation} 
\end{array} \right)  \nonumber \\
= \left(
\begin{array}{c | c}
\Yabc &  \Yin \\
\hline 
\Yin &  \Ynn 
\end{array}
\right) .
\end{IEEEeqnarray} 

Similarly, the sequence admittance matrix $\Yseq$ can be obtained by performing symmetric component transformation over the partitioned shunt admittance $\Yabcn$ \cite{Fred_transformation_2022},
\begin{IEEEeqnarray}{c}
\Yseq = \Ainv \Yabc \A .
\end{IEEEeqnarray}

\subsection{Expressions for $\GMR$ and $\Dij$ w.r.t. Coordinates}
Carson's equations require values of $\GMR$ and $\Dij$ to derive the impedance. 
The $\GMR$ of a conductor is related to its strand radius $\radius$ by a ratio $\Kgmr$. 
The ratio $\Kgmr$ depends on the number of strands $\N$ and strand geometry \cite{Zubair2014ArithmeticAL}, 
\begin{IEEEeqnarray}{c}
\GMR =\Kgmr \cdot \radius ,  \quad \quad \,\, \radius \,\&\, \GMR \geq 0 \label{GMR eq}.
\end{IEEEeqnarray}
We assume all conductors have the same strand radius $\radius$ and hence $\GMR$. To compute $\Dij$, Cartesian coordinates are adopted and $\Dij$ is defined to be the distance between the center of circular stranded conductors. Although this approach is a simplified alternative to calculate the geometric mean distance (GMD) between strands of conductors, Urquhart and Thomson report that this approach only gives negligible difference compared with GMD and is less subject to rounding error \cite{urquhart_series_2015}.
For a system with $\n$ conductors, the x- and y-coordinate of circular-stranded conductors are contained in vectors $\x$ and $\y$ of length $\n$, respectively,
\begin{align}
\x &=  \begin{pmatrix}
x_{\config,1} &
\hdots & 
x_{\config,n} 
 \end{pmatrix}^{\text{T}},
 &
 \y &= \begin{pmatrix}
y_{\config,1} &
\hdots &
y_{\config,n} 
\end{pmatrix}^{\text{T}} .\label{eq_xy_definition}
 \end{align}
The (natural) Euclidean and mirror-image `distances' are,
\begin{IEEEeqnarray}{c}
 \Dij = \sqrt{(x_{\config,i} - x_{\config,j})^2 + (y_{\config,i} - y_{\config,j})^2 },  \label{Euclidean equation}  \IEEEyesnumber \IEEEyessubnumber\\
 \Sij = \sqrt{(x_{\config,i} - x_{\config,j})^2 + (y_{\config,i} + y_{\config,j})^2 } \label{image Euclidean equation}.  \IEEEyessubnumber
\end{IEEEeqnarray}

\subsection{\textcolor{\varcolor}{Likely Geometries for Overhead Lines}} \label{section: OH geometry}
\begin{figure*}[th!]
\centering
\includegraphics[width=\linewidth]{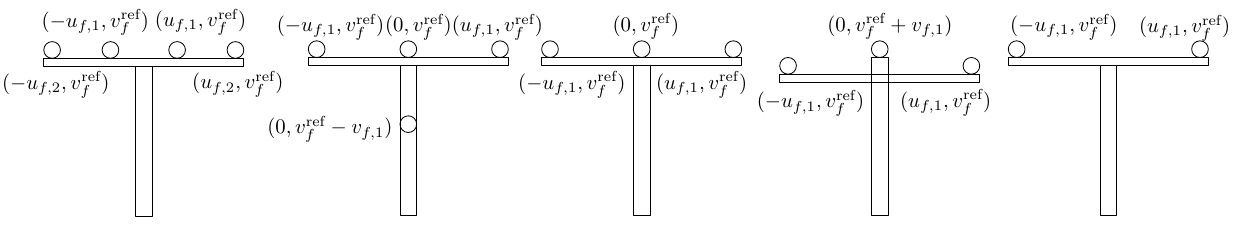}
\caption{Overhead line geometries; from left to right:  4-wire horizontal plane, 4-wire neutral-under, 3-wire flat intermediate, 3-wire triangular arrangement and 2-wire horizontal plane.}
\label{fig:OH}
\end{figure*}

For overhead lines, there are many pole geometries, but in this work they are limited to the most common configurations depicted in Fig. \ref{fig:OH}. 
For the sake of convenience, free variables $u_{\config,1}\leq u_{\config,2} \leq \dots \leq u_{\config,n}$ are taken to describe horizontal coordinates in $\x$ and similarly $v_{\config,1} \leq v_{\config,2} \leq \dots \leq v_{\config,n}$ for vertical coordinates in $\y$, with all of them being nonnegative. 

Next,  geometries follow some practical constraints. Firstly, air is assumed to be the only insulation medium so a minimal distance $\Dmin$ between wire applies. Next, the crossarm has a finite length, which is an upper bound for the distance between wires, $\umaxOH$. Lastly, all  wires share the same reference distance above ground $\vref$, which is subject to height limitation,  
\begin{IEEEeqnarray}{c}
\vrefmin \leq \vref \leq \vrefmax. \label{eq_vref_bound}
\end{IEEEeqnarray}
\subsubsection{Horizontal Plane (4-wire)} This is the simplest geometry of four wires such that wires are placed apart from the center of the pole by either $\freevar{u}{\config,1}$ or $\freevar{u}{\config,2}$. 
Neutral wire is assumed to be located furthermost from the pole center, e.g. $(x_4,y_4)=(\freevar{u}{\config,2},\vref)$, and the coordinate vectors are,   
\begin{IEEEeqnarray}{c}
 \x
 =
 \begin{pmatrix}
-\freevar{u}{\config,2} &
-\freevar{u}{\config,1} &
\freevar{u}{\config,1} &
\freevar{u}{\config,2} 
 \end{pmatrix}^\text{T}, \label{eq_x_horizontal}
 \IEEEyesnumber \IEEEyessubnumber\\
  \y
 = \vref \cdot
 \begin{pmatrix}
1 &
1 &
1 &
1 
 \end{pmatrix}^\text{T} \label{eq_oh_horizontal}. \IEEEyessubnumber
\end{IEEEeqnarray}
The minimum distance constraint is defined  on the horizontal coordinates. The two outermost wires must be  within the half-side crossarm length of $\umaxOH$, and separated from the two wires near the center of the pole with at least $\Dmin$,  
\begin{IEEEeqnarray}{c}
\umaxOH \geq \freevar{u}{\config,2} 
 \geq \freevar{u}{\config,1} + \Dmin \label{eq_4OH_hor_cons_u2} .
 \end{IEEEeqnarray}
Next, the two wires near the center of the pole must be separated from each other by at least $\Dmin$,
\begin{IEEEeqnarray}{c}
\freevar{u}{\config,1} \geq \Dmin/2.  \label{eq_4OH_hor_cons_u1}
\end{IEEEeqnarray}

\subsubsection{Neutral-under (4-wire)} 
This geometry has three phase wires on the top and a neutral wire under the phase wires, 
\begin{IEEEeqnarray}{c}
 \x
 =
 \begin{pmatrix}
-\freevar{u}{\config,1} &
0 &
0 &
\freevar{u}{\config,1} 
 \end{pmatrix}^\text{T}, \IEEEyesnumber \IEEEyessubnumber\\
    \y 
 = 
 \begin{pmatrix}
\vref &
\vref &
\vref-\freevar{v}{\config,1} &
\vref
 \end{pmatrix}^\text{T} \label{eq_neutral_under}. \IEEEyessubnumber
\end{IEEEeqnarray}
Distances are upper bounded by crossarm length and maintain minimal distance of $\Dmin$ from each other,
\begin{IEEEeqnarray}{c}
\umaxOH \geq \freevar{u}{\config,1} \geq \Dmin. \label{eq_4OH_neutral_under_cons_u1} 
\end{IEEEeqnarray}
The neutral wire is above ground and  below the center phase wire by at least $\Dmin$,  
\begin{IEEEeqnarray}{c}
\vref \geq \freevar{v}{\config,1} \geq \Dmin \label{eq_4OH_neutral_under_cons_vref}.
\end{IEEEeqnarray}

\subsubsection{Horizontal Plane (3-wire)} 
One wire is placed at the top center of the pole and the remaining two away from the center with distance $\freevar{u}{\config,1}$ ,
\begin{align}
 \x
 &=
 \begin{pmatrix}
-\freevar{u}{\config,1}& 
0 &
\freevar{u}{\config,1} 
 \end{pmatrix}^\text{T},
 &
  \y
 &= \vref 
 \begin{pmatrix}
1 &
1 &
1 
 \end{pmatrix}^\text{T}.
\end{align}
In this configuration, distance constraint \eqref{eq_4OH_neutral_under_cons_u1} also applies. 

\subsubsection{Triangular Arrangement (3-wire)} 
In this geometry three wires form an isosceles triangle such that one wire is placed at the center of the pole above the crossarm, and two wires are placed $\freevar{u}{\config,1}$ away from the center of pole, 
\begin{align}
 \x
 &=
 \begin{pmatrix}
-\freevar{u}{\config,1} &
0 &
\freevar{u}{\config,1} 
 \end{pmatrix}^\text{T} ,
 &
  \y
 &= 
 \begin{pmatrix}
\vref &
\vref+\freevar{v}{\config,1} &
\vref \\
 \end{pmatrix}^\text{T} ,
\end{align}
\textcolor{\varcolor}{
where $\tan(\theta)$ is the ratio between $\freevar{u}{\config,1}$ and $\freevar{v}{\config,1}$},
\begin{IEEEeqnarray} {c}
\freevar{v}{\config,1}=\freevar{u}{\config,1} \cdot \tan (\theta) \label{eq_3OH_tri_cons_theta} .
\end{IEEEeqnarray}
The two wires on the crossarm are separated by at least $\Dmin$ and each wire on the crossarm should be at least $\Dmin$ away from the wire on the top of the pole. Finally, the lower bound on $\freevar{u}{\config,1}$ needs to consider both conditions,
\begin{IEEEeqnarray} {c}
\umaxOH \geq \freevar{u}{\config,1} \geq \max \{ \Dmin/2,\Dmin \cdot \cos {\theta}   \}    \label{eq_3OH_tri_cons_max}.
\end{IEEEeqnarray}

\subsubsection{Horizontal Plane (2-wire)} This is a common arrangement for single-phase connection. Two lines are equally separated from the middle of the pole,  
\begin{align}
 \x
 &=
 \begin{pmatrix}
-\freevar{u}{\config,1} &
\freevar{u}{\config,1} 
 \end{pmatrix}^\text{T} ,
 &
  \y
 &= 
\vref  \begin{pmatrix}
1 &
1
 \end{pmatrix}^\text{T} .
\end{align}
Similarly, distance bound apply for $\freevar{u}{\config,1}$, 
\begin{IEEEeqnarray}{c}
\umaxOH \geq \freevar{u}{\config,1} \geq \Dmin/2 .
\end{IEEEeqnarray}

\subsection{\textcolor{\varcolor}{Likely Geometries for Low-voltage Aerial Bundle Cable (LVABC) and Underground Cable (UGC) }}\label{section: cable geometry}
\begin{figure}[tbh]
\centering
\includegraphics[width=1\linewidth]{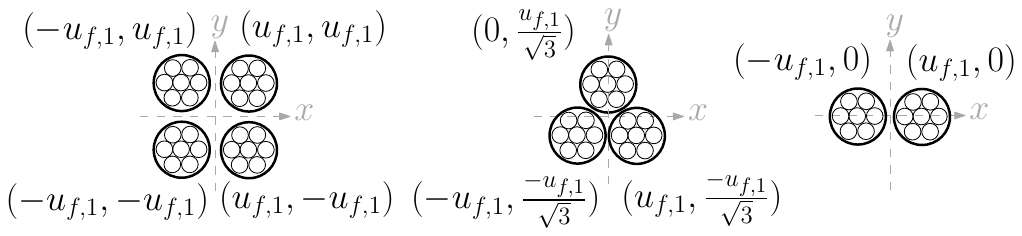}
\caption{Cable geometry normalised in $\freevar{u}{\config,1}$: from left to right are 4-core square, 3-core equilateral triangle and 2-core closely-packed horizontal. Reference height of cables are not shown for simplicity of illustration.}
\label{fig:cable}
\end{figure}

Common geometries of 4-core, 3-core and 2-core cables are shown in Fig. \ref{fig:cable}. 
It is assumed that radius $\Rnom$ of all cores are the same and all cores are closely packed symmetrically in a cable. 
Common cables type in low-voltage application are circular stranded with each core insulated with thickness $\tnom$. 
Because all cores in a cable are closely-packed, it is convenient to express the Cartesian coordinates in terms of core radius $\Rnom$. To begin with, $\Rnom$ is derived from the strand radius $\radius$ and the core insulation thickness $\tnom$,
\begin{IEEEeqnarray}{c}
\Rnom=\R + \tnom, \quad \R=\Kradius \cdot \radius \label{eq_Rnom}.
\end{IEEEeqnarray}

\textcolor{\varcolor}{The packing coefficient $\Kradius$ is the coefficient relating strand radius and radius of core excluding the core insulation thickness. The term $\Kradius \cdot \radius$ is the minimal radius of core required to pack all strands without insulation.}

To be consistent with overhead notation for distance, $\Rnom$ is set to equal to $\freevar{u}{\config,1}$,
\begin{IEEEeqnarray}{c}
\freevar{u}{\config,1}=\Rnom \label{eq_Rnom_u1}.
\end{IEEEeqnarray}
\textcolor{\varcolor}{
Lastly, $\ua$ is bounded over some practical values,
\begin{IEEEeqnarray}{c}
\umincable \leq \ua \leq \umaxcable.
\end{IEEEeqnarray}
}
\subsubsection{4-conductor} This is the geometry for both 4-core circular- and sector-stranded cable\footnote{\textcolor{\varcolor}{Sector-stranded geometry follows Fig. 9b of \cite{Zubair2014ArithmeticAL}, with $\GMR$ and $\Dij$ calculated using (10) and (14) of \cite{cleenwerck_approach_2022}. This results in slightly rectangular geometry in $\Dij$. Coefficient $\Kradius$ for deriving $\R$ for shunt computation is left for future work. In inverse estimation, we assume 4-core geometry for simplification of analysis. Eq. \eqref{eq_Rnom} does not hold for sector-stranded cable.}}, with reference height $\vref$ added to the local y-coordinate, 
\begin{IEEEeqnarray}{c}
 \x
 = \freevar{u}{\config,1} \cdot
 \begin{pmatrix}
1 &
-1 &
-1 &
1 
 \end{pmatrix}^\text{T} , \IEEEyesnumber \IEEEyessubnumber\\
   \y
 = \freevar{u}{\config,1}  \cdot
 \begin{pmatrix}
1 &
1 &
-1 &
-1 
 \end{pmatrix}^\text{T} \IEEEyessubnumber
+\vref.
\end{IEEEeqnarray}

\subsubsection{3-conductor} 
The 3 cores form an equilateral triangle,
\begin{IEEEeqnarray}{c}
 \x
 =\freevar{u}{\config,1}  \cdot
 \begin{pmatrix}
-1 &
0 &
1 
 \end{pmatrix}^\text{T} ,  \IEEEyesnumber \IEEEyessubnumber\\
  \y
 ={\freevar{u}{\config,1}}/{\sqrt{3}}  \cdot
 \begin{pmatrix}
-1 &
2 &
-1 
 \end{pmatrix}^\text{T} \IEEEyessubnumber 
 +\vref.
 \label{eq_3_cable_cable}
\end{IEEEeqnarray}

\subsubsection{2-conductor} The 2 cores are adjacent, 
\begin{IEEEeqnarray}{c}
 \x
 =\freevar{u}{\config,1} \cdot
 \begin{pmatrix}
-1 &
1 
 \end{pmatrix}^\text{T}, \IEEEyesnumber \IEEEyessubnumber\\
  \y
 = \vref \cdot 
 \begin{pmatrix}
1 &
1 
 \end{pmatrix}^\text{T}. \IEEEyessubnumber
\end{IEEEeqnarray}

\subsection{Electrical Characteristic Modelling of Material} \label{sec:resistivity}
The IEC 60287-1-1 standard is adopted for deriving the ac resistance and details of function $h$ in \eqref{eq_rac_dummy} is now expanded.  
For conductors with single material $\material$, e.g. Aluminium or Copper, the ac resistance $\rac$ is derived as follows. Firstly, the cross sectional area $\area$ of circular-stranded conductor is the sum of cross sectional area of $\N$ strands,
\begin{IEEEeqnarray}{c}
\area = \N \pi \radius^2 \label{eq_area} .
\end{IEEEeqnarray}
Then, the dc resistance of conductor $\rdc$ at 20\,$^\circ$C is obtained through dividing the material resistivity $\p$ at 20\,$^\circ$C by the cross sectional area $\area$.
Next, $\rdc$ at the expected temperature $\T$ is derived using  the material temperature coefficient $\Tcoef$,
\begin{IEEEeqnarray}{c}
\rdc = \frac{\p}{\area} \left(1+\Tcoef \left(\T-20^\circ \text{C} \right) \right)  \label{eq_rdc}.
\end{IEEEeqnarray}
We assume all conductors share the same temperature $\T$. Now, ac resistance $\rac$ is calculated by considering increase in resistance due to skin effect $\Cskin$ and proximity effect $\Cprox$,
\begin{IEEEeqnarray}{c}
\rac=(1+\Cskin)(1+\Cprox)\rdc \label{eq_rac} ,
\end{IEEEeqnarray}
\textcolor{\varcolor}{The details for deriving $\Cskin$ and $\Cprox$ are in IEC60287-1-1 \cite{IEC60287}.}

The range of radius, area, temperature and \textcolor{\varcolor}{insulation thickness}  of conductors are bounded over practical values for the inverse problem,
\begin{IEEEeqnarray}{c}
\rmin \leq \radius \leq \rmax  \IEEEyesnumber \IEEEyessubnumber,\label{eq_r_bound}\\
\Amin \leq \area \leq \Amax ,\IEEEyessubnumber \label{eq_A_bound}\\
\Tmin \leq \T \leq \Tmax , \IEEEyessubnumber \label{eq_T_bound} \\
\textcolor{\varcolor}{\tnommin \leq \tnom \leq \tnommax .}\IEEEyessubnumber \label{eq_tnom_bound}
\end{IEEEeqnarray}

\section{Inverse Problem as an Optimisation Problem}\label{sec_opt_prob_spec}
\textcolor{\varcolor}{We assume utilities only provide reference sequence components $\Rzeroref, \Roneref, \Xzeroref, \Xoneref, \Bzeroref, \Boneref$ as derived from the unbalanced line\footnote{Expection is balanced 3-phase cable with equilateral triangle geometry.} using the mathematical steps from Section \ref{section: network modelling}. No information of discrete properties material $\material$ and configuration $\config$ (which contains discrete information of number of conductors $n$, number of strands $N$ and type of geometry) are provided. Therefore, they need to be considered individually such that all discrete properties are combined to form combination $ \combination \in \setcombination \subseteq \setlineid \times \setconfig \times  \setmaterial$, which is optimised to match the reference sequence components.}

As illustrated in Fig. \ref{fig:flow_block_diagram}, the inverse problem is split into three sequential steps, each cast as optimisation problems, namely 1) the feasibility problem, 2) the bound tightening problem and 3) the bounded slack problem. In the feasibility problem, a combination $\combination$ is optimised to return the sequence components closest to the reference sequence components given. The closer it is, the more likely a combination is a correct estimation. In the bound tightening problem, the range of decision variables returned from from feasibility problem is studied. This provides numerical evidence of the uniqueness of the solution to the feasibility problem. Lastly, in the bounded slack problem, sensitivity of decision variables over reference sequence components will be analysed. It studies the range of decision variables when there is slack in the reference sequence components. \textcolor{\varcolor}{Also, it investigates how much slack is required over the reference sequence components so that an incorrect combination becomes feasible to match the reference sequence components.} It provides an understanding how resilient decision variables are when potential error are introduced in sequence component data, such as using wrong fundamental frequency (50 or 60\,Hz) or different variants of Carson's equations.    

\subsection{Feasibility Problem}
The distance of diagonal impedance and susceptance sequence components $\BZdiff$ between a combination $\combination$ and the reference components is minimised. 
The error terms within the absolute operators in \eqref{obj_eq_feasibility_BZ_diff} represent $\ell_1$ norms, which induce sparsity in the error terms, to minimise error propagation across different sequence component values. Next, it is normalised by the corresponding reference sequence component value and divided by 6, i.e.,
\begin{IEEEeqnarray}{c}
\BZdiff = \Big( \frac{1}{ \Bzeroref} \underbrace {|\Bzero-\Bzeroref|}_{\Bzeroaux} + 
\frac{1}{\Boneref} \underbrace {|\Bone-\Boneref|}_{\Boneaux} \nonumber + \\
\frac{1}{ \Rzeroref} \underbrace {|\Rzero-\Rzeroref|}_{\Rzeroaux} +
\frac{1}{\Roneref} \underbrace {|\Rone-\Roneref|}_{\Roneaux} \nonumber + \\
\frac{1}{\Xzeroref} \underbrace{|\Xzero-\Xzeroref|}_{\Xzeroaux} +
\frac{1}{\Xoneref}  \underbrace{|\Xone-\Xoneref|}_{\Xoneaux}\Big)/6.\label{obj_eq_feasibility_BZ_diff}
\end{IEEEeqnarray}
In case of missing shunt susceptance data, the corresponding terms are dropped from \eqref{obj_eq_feasibility_BZ_diff}, and the division by 6 is replaced with division by 4 (now labelled $\Zdiff$).

\begin{IEEEeqnarray}{l c}
\min & \BZdiff,  \label{eq_feasibility}\\
\text{var. }  & \x,  \y, \radius, \Dij,  \Sij, \T, \rdc, \rac, \area, \R, \textcolor{\varcolor}{\Rnom,} \nonumber \\
 & \textcolor{\varcolor}{\vref,\tnom,} \GMR,\Rabcn, \Xabcn, \Rabc, \Xabc, \Rseq, \Xseq,  \nonumber  \\ 
\text{s.t. } & \eqref{eq_Z_abcn} \text{ to } \eqref{obj_eq_feasibility_BZ_diff}. \nonumber
\end{IEEEeqnarray}
The optimal values of $R_{l,00}, X_{l,00}, R_{l,11}, X_{l,11}, B_{l,11}, B_{l,11}$, as denoted by $\Rzerostar,\Ronestar,\Xzerostar,\Xonestar,\Bzerostar,\Bonestar$ respectively, will be passed to bound-tightening problem.

\subsection{Bound Tightening Problem}
The range of decision variables is investigated by evaluating their minimal and maximal values. Each decision variable $\var \in \setvar=\{\radius,\T,\freevar{u}{\config,1} ,\freevar{u}{\config,2} ,\freevar{v}{\config,1}, \textcolor{\varcolor}{\vref\}}$ is iterated individually. A negative sign in \eqref{eq_obj_bound_tightening} effectively indicates maximisation. Note that the optimal sequence components values with asterisk obtained from the previous feasibility problem will be used here as shown in \eqref{eq_BZ_asteriska} - \eqref{eq_BZ_asterisk}.
\begin{IEEEeqnarray}{l c}
\underset{\var \in \setvar}{\min} & \pm \text{\,} \var,     \label{eq_obj_bound_tightening} \IEEEyesnumber  \IEEEyessubnumber\\ 
\text{var. }  & \x,  \y, \radius, \Dij,  \Sij, \T, \rdc, \rac, \area, \R, \textcolor{\varcolor}{\Rnom,} \nonumber \\
 & \textcolor{\varcolor}{\vref,\tnom,} \GMR,\Rabcn, \Xabcn, \Rabc, \Xabc, \Rseq, \Xseq,  \nonumber  \\ 
\text{s.t. } & \eqref{eq_Z_abcn} \text{ to } \eqref{eq_tnom_bound}, \nonumber\\
&\Rzero = \Rzerostar,
\Rone = \Ronestar, \IEEEyessubnumber \label{eq_BZ_asteriska} \\
&\Xzero = \Xzerostar,
\Xone = \Xonestar,  \IEEEyessubnumber \\
&\Bzero = \Bzerostar,
\Bone = \Bonestar.  \IEEEyessubnumber  \label{eq_BZ_asterisk}
\end{IEEEeqnarray}

Note that only diagonal entries are fixed while the off-diagonal entries in $\Rseq, \Xseq$ remain free.

\subsection{Bounded Slack Problem}
Sensitivity analysis of decision variables $\var$ in $\setvar$ over reference sequence components are conducted by introduction of slack $\slack$ from shown in \eqref{eq_slack_R0} to \eqref{eq_slack_B1}.
\begin{IEEEeqnarray}{l c}
\underset{\var \in \setvar}{\min} & \pm \text{\,} \var,  \label{eq_obj_bound_tightening_slack} \IEEEyesnumber  \IEEEyessubnumber\\ 
\text{var. }  & \x,  \y, \radius, \Dij,  \Sij, \T, \rdc, \rac, \area, \R, \textcolor{\varcolor}{\Rnom,} \nonumber \\
 & \textcolor{\varcolor}{\vref,\tnom,} \GMR,\Rabcn, \Xabcn, \Rabc, \Xabc, \Rseq, \Xseq,  \nonumber  \\ 
\text{s.t. } & \eqref{eq_Z_abcn} \text{ to } \eqref{eq_tnom_bound}, \nonumber \\
& (1+\slack)\Rzeroref \geq \Rzero \geq (1-\slack)\Rzeroref,  \IEEEyessubnumber \label{eq_slack_R0}\\
& (1+\slack)\Roneref \geq \Rone \geq (1-\slack)\Roneref,  \IEEEyessubnumber \label{eq_slack_R1}\\
& (1+\slack)\Xzeroref \geq \Xzero \geq (1-\slack)\Xzeroref,  \IEEEyessubnumber \label{eq_slack_X0}\\
& (1+\slack)\Xoneref \geq \Xone \geq (1-\slack)\Xoneref,  \IEEEyessubnumber \label{eq_slack_X1}\\
& (1+\slack)\Bzeroref \geq \Bzero \geq (1-\slack)\Bzeroref,  \IEEEyessubnumber \label{eq_slack_B0}\\
& (1+\slack)\Boneref \geq \Bone \geq (1-\slack)\Boneref.  \IEEEyessubnumber \label{eq_slack_B1}
\end{IEEEeqnarray}

Slack  $\beta$ is introduced in the reference sequence components to study when incorrect inverse estimation becomes feasible. 
\subsection{Reformulation as Differentiable Equations in the Reals}
Optimization modeling toolboxes (typically) do not support complex numbers, so for implementation purposes equations are re-stated in real variables. Next, the square root and absolute value functions are non-differentiable whereas division operation requires addition of nonlinear constraint, so they are reformulated as quadratic equations.

Firstly, the absolute value terms is reformulated by introducing auxiliary variables to represent the value of the absolute value of the difference, i.e.,
\begin{IEEEeqnarray}{c}
\min |a-b|, \text{s.t. } a,b \in \mathbb{R}, \nonumber \\
\iff \min c, \text{s.t. }  c \geq a-b, c\geq b-a  \implies c = |a-b| \nonumber .
\end{IEEEeqnarray}

Next,  Kron's reduction \eqref{Kron reduction eq} and the symmetrical component transform \eqref{SCT eq} are equivalently restated in the reals as quadratic  and linear equations respectively. 


The square-root is removed by squaring both side of \eqref{Euclidean equation} and \eqref{image Euclidean equation},
\begin{IEEEeqnarray}{c}
 (\Dij)^2 = (x_{\config,i} - x_{\config,j})^2 + (y_{\config,i} - y_{\config,j})^2, \Dij\geq 0, \label{eq_dist_differentiable}\\
 (\Sij)^2 = (x_{\config,i} - x_{\config,j})^2 + (y_{\config,i} + y_{\config,j})^2 , \Sij\geq 0   \label{eq_dist_mirror_differentiable} .
\end{IEEEeqnarray}
Lastly, division operation in \eqref{eq_rdc}  is removed by rearrangement of variables,
\begin{IEEEeqnarray}{c}
\rdc \cdot \area = \rho_\material (1+\Tcoef (\T-20)),  \,\,\,\,\,\, \rdc \geq 0 \label{eq_rdc_new}.
\end{IEEEeqnarray}

\section{Test case and validation}\label{section: test_case_and_validation}
\textcolor{\varcolor}{
We perform numerical studies using a dataset of  common OH line construction types and cables based on \cite{Ergon_overhead,Ergon_underground}. 
Firstly, we conduct simple test case by generating $\Rzeroref$, $\Roneref$, $\Xzeroref$, $\Xoneref$, $\Bzeroref$ and $\Boneref$ for some of the entries in our dataset.
Next, the recovery algorithm is applied to provide insight on what information can be reliably recovered and the risk of mismatch. 
Then, across the whole dataset we validate insights from the initial examples. Lastly, recovery over utility data is performed to demonstrate a real-world application.}

The inverse problem is set up in Julia (v1.8.3) with optimisation modelling language JuMP and solver Ipopt. Ipopt in combination with JuMP's automatic differentiation allows the formulation of optimisation problems with quadratic and transcendental constraints. The computation is performed on an Intel(R) Core(TM) i7-11700 @ 2.50\,GHz  with  32\,GB RAM. The forward calculation result is validated against OpenDSS with results matched up to 4 decimal places in [$\Omega$/km] for impedance and \textcolor{\varcolor}{2 decimal places in [$\mu S$/km] for admittance.}

\subsection{Data on Overhead Lines and Cables}
The test case as summarised in Table \ref{tab_EQLD_construction_code} considers overhead line and cables commonly adopted in LV network in Eastern Australia \cite{Ergon_overhead,Ergon_underground}, which the authors of this paper have knowledge on.

\begin{table}[tbh]
    \caption{Common overhead line and cable for LV network in Queensland Australia \cite{Ergon_overhead}\cite{Ergon_underground}.}
    \label{tab_EQLD_construction_code}
    \centering
    \begin{tabular}{l l l l l l l}
        \toprule
        $\lineid \in \setlineid $ & $|\mathcal{W}|$ & $\area$ &  \textcolor{\varcolor}{$\rstd$} & $\N$ &  $\material$ & $\tnom$ \\ 
        & & [mm$^2$] & [mm] & & & [mm]
        \\
        \midrule
        OH Libra & 3 or 4 & 49.48  & 1.5 & 7 & Al-1350 & N/A \\
        OH Mars & 3 or 4 & 77.31  &  1.875 & 7 & Al-1350 & N/A \\
        OH Moon & 3 or 4 & 124.04  &  2.375 & 7 & Al-1350 & N/A \\
        LVABC4x95 & 4 &  94.75 & 1.26 & 19 & Al-1350 & 1.7 \\  
        LVABC4x50  & 4 &  48.17 & 1.48 & 7 & Al-1350 & 1.5 \\
        LVABC4x25  & 4 &  26.61& 1.1 & 7 & Al-1350 & 1.3 \\
        LVABC3x25  & 3 &  26.61 & 1.1 & 7 & Al-1350 & 1.3 \\
        UGC16x4Cu & 4 & 15.89 & 0.85 & 7 & Cu & 1 \\
        UGC50x4Cu & 4 & 48.17  & 1.48 & 7 & Cu & 1.5 \\
        UGC240x4Al & 4 &  239.40 & 1.26 & 48 & Al-1350 & 1.7 \\
        \bottomrule
    \end{tabular}
\end{table} 

Two-wire overhead lines and two-core cables are excluded because using sequence coordinates to represent their impedance is misleading. Transformer impedance is out-of-scope. Next, skin and proximity effect are excluded by setting $\Cskin=\Cprox=0$ to avoid the associated non-linear equations as per IEC 60287-1-1 standard\cite{IEC60287}. This is a mild assumption based on the study from Urquhart and Thomson that for a large 300 mm$^2$ cable, the combined skin and proximity effect together account for only 2$\%$ in ac resistance\cite{urquhart_series_2015}. \textcolor{\varcolor}{Bounds for variables are added as shown in Table \ref{tab_bound_value_code}. Lastly, superscript ``std'' is added to variables $\radius,\ua,\ub,\va,\vref$, which corresponds to standardised strand radius and geometry parameters as adopted by manufacturers or utilities.}

\begin{table}[tbh]
    \setlength{\tabcolsep}{1.2pt}
    \caption{Bound for optimisation variables.}
    \label{tab_bound_value_code}
    \centering
    \begin{tabular}{l l l}
    \toprule
    Parameter &  Value & Justification \\
    \midrule
    $\Dmin$ & 380\,mm & Ausgrid standard \cite{ns220_nodate}.\\
    $\umaxOH$ & 1500\,mm & Energy QLD largest crossarm \cite{Ergon_overhead}.\\ 
    $\umincable$ & 2.55\,mm & UGC16x4Cu core radius w/o $\tnom$. \cite{Nexan_cable_catalogue}.\\ 
    $\umaxcable$ & 30\,mm & UGC300x4Al cable radius in \cite{sector_catalogue_NZ}. \\
    $\Tmin$ & 0\,$^\circ$C & Arbitiary for simulation purpose.\\
    $\Tmax$ & 105\,$^\circ$C & Emergency oper. temp. for cable \cite{Ergon_underground}.\\
    $\rmin,\rmax$ & 0.85, 2.375\,mm & Min. and max. $\radius$ in Table \ref{tab_EQLD_construction_code}. \\
    \textcolor{\varcolor}{$\tnommin,\tnommax$} & \textcolor{\varcolor}{1, 1.7\,mm} & \textcolor{\varcolor}{Min. and max. $\tnom$ in Table \ref{tab_EQLD_construction_code}.} \\
    $\Amin$,$\Amax$ &15, 240\,mm$^2$ & Min. and max. area in Table \ref{tab_EQLD_construction_code}.\\
    $\Amin$,$\Amax$ &185, 300\,mm$^2$ & Sector area of $\N$=37 \& $\N$=61\cite{sector_catalogue_NZ}. \\
    \textcolor{\varcolormoving}{$\vrefmin,\vrefmax$} & \textcolor{\varcolormoving}{5.8, 21.5\,m} & \textcolor{\varcolormoving}{Energy QLD OH line standard \cite{Ergon_overhead}.} \\
    \textcolor{\varcolormoving}{$\vrefmin,\vrefmax$} & \textcolor{\varcolormoving}{-6, -0.6\,m} & \textcolor{\varcolormoving}{Energy QLD UGC standard \cite{Ergon_underground}.} \\
    \bottomrule
    \end{tabular}
\end{table}

\subsection{Scenario Setting}
From Table \ref{tab_EQLD_construction_code}, overhead lines in Australia are typically 7-strand Aluminium wire. Two 4-wire geometries and three 3-wire geometries, as summarised on the left of Table \ref{tab_forward_OH}\footnote{\textcolor{\varcolormoving}{For neutral-under in Table \ref{tab_forward_OH}, $\uastd$ and $\vastd$ are the averaged values of Fig. 5 of \cite{kersting_distribution_1995} and its mirror image.}}, are adopted in Australia. Therefore, overhead lines only differ in configuration $\config$. For cables, it has different number of cores (3 / 4), number of strands (7 / 19 / 48) and material (Al-1350 / Cu) in Table \ref{tab_material_parameters}, therefore cables differ in both configuration $\config$ (which contains $\n$ and $\N$) and material $\material$.
In inverse estimation, assuming there is no knowledge of which configuration are the sequence components generated from (and no knowledge of material if it is a cable) but knowing whether it is a cable or overhead line. Combinations are therefore formed by iterating the reference sequence components by configurations in Table \ref{tab_config_parameters} (and material types if it is a cable). 
We assume utility to provide all diagonal sequence components and know whether they belong to an overhead line or cable.
\begin{table}[tbh]
    \textcolor{\varcolor}{
    \setlength{\tabcolsep}{3pt}
    \caption{Standardised OH geometries \cite{kersting_distribution_1995,dugan_open_2011} with $\vrefstd$=9150[mm] \cite{Ergon_overhead}. On the right is Forward calculation with Mars at 75\,°C.}
    \label{tab_forward_OH}
    \centering
    \begin{tabular}{l l l l l |l}
    \toprule
    Geometry & $|\mathcal{W}| $ & $\uastd$ & $\ubstd$ & $\vastd$ & ($\Rzeroref,\Xzeroref,\Roneref,\Xoneref$)\\
             &     &   [mm]       &    [mm]      &    [mm]      & [$\Omega/km$] \\
    \midrule
    Hori. (4w)        & 4	& 450  & 1100&  N/A  & (0.7788,1.1057,0.4481,0.3422) \\ 
    Neutral-under  & 4	& 1118 & N/A &  1575 & (0.7554,1.1072,0.4472,0.3671) \\ 
    Hori. (3w)	       & 3	& 1100 & N/A &  N/A	& (0.5952,1.5934,0.4472,0.3662) \\ 
    Tri. ($\theta$=21.67°)  & 3	& 1100 & N/A &  N/A	& (0.5952,1.5873,0.4472,0.3692) \\ 
    Tri. ($\theta$=49.27°)  & 3	& 508  & N/A &  N/A	& (0.5952,1.6547,0.4472,0.3355) \\ 
    \bottomrule
    \end{tabular}
    }
\end{table}

\begin{table}[tbh]
\textcolor{\varcolor}{
\setlength{\tabcolsep}{2.5pt}
\caption{material and configuration parameters.}
\label{tab_material_parameters}
    \centering 
    \begin{tabular}{l l l l |l l l}
        \toprule
        $m \in \mathcal{M}$ & $\p$& $\Tcoef$ & Standard & $\N\in\setconfig$ & $\Kgmr$ & $\Kradius$ \\
        & [$10^{-9} \Omega\,m$] &[$1/^{\circ}C$] & & & &\\
        \midrule
        Al-1350 & 28.3	& 0.00403 & AS3607-1989\cite{AS3607} & 7 & 2.18 & 3\\
        Cu & 17.77 & 0.00381 & AS1746-1991\cite{AS1746} & 19 & 3.79 & 5\\
         &  &  &                                & 48 & 6.41 & [-]\\
        \bottomrule
    \end{tabular}
}
\end{table}

\begin{table}[tbh]
\setlength{\tabcolsep}{7pt}
\textcolor{\varcolor}{
\caption{Configurations for inverse estimation.}
\label{tab_config_parameters}
    \centering
    \begin{tabular}{l l l |l l l }
        \toprule
$\config \in \setconfig^\text{OH}$ &  $|\mathcal{W}| $ & $\N$ & 
$\config \in \setconfig^\text{cable}$ &  $|\mathcal{W}| $ & $\N$ \\
        \midrule
        Hori. (4w)             & 4 & 7  &  3w7N circ.  & 3 & 7 \\
        Neutral-under          & 4 & 7  &  3w19N circ.  & 3 & 19 \\
        Hori. (3w)             & 3 & 7   &  4w7N circ.  & 4 & 7 \\
        Tri. ($\theta$=21.67°) & 3 & 7  &  4w19N circ.  & 4 & 19 \\
        Tri. ($\theta$=49.27°) & 3 & 7  &  4w48N sect.  & 4 & 48 \\ 
        \bottomrule
    \end{tabular}
    }
\end{table}

\textcolor{\varcolor}{
\subsection{Simple Test Cases}
For overhead lines, as they differ in geometry types, we firstly perform forward calculation to understand potential mismatch due to similar diagonal sequence components, as shown on the right of Table \ref{tab_forward_OH}. Large differences exist in $\Rzeroref$ and $\Xzeroref$ between 3-wire and 4-wire geometries, suggesting that number of wires $\n$ is likely to be reliably distinguished. For standardised geometries with the same number of wires, in general only minor differences are observed for sequence reactance, especially between Hori. (3w) and Tri. ($\theta$=21.67°). However, the risk of mismatch is low because of the small difference in reactance. In contrast, standardised Tri. ($\theta$=49.27°) shows a considerably different $\Rzeroref$ and $\Roneref$ over the other two 3-wire geometries, because of its much smaller $\uastd$ leading to considerably smaller conductor spacing distances. This allows us to distinguish it from the other two.}

\textcolor{\varcolor}{
Now the recovery algorithm is applied to standardised Tri. ($\theta$=21.67°) with conductor Mars. The result is summarised in Table \ref{tab_inverse_OH}. In the feasibility problem, $\Zdiff$ is 0 for all 3-wire geometries, but are 0.137 and 0.0653 for Hori. (4w) and neutral-under geometry respectively. The 4-wire combinations are therefore reliably eliminated. In the bound-tightening problem, local uniqueness of $\ua$ and $\radius$ are observed and correctly recovered for correct inverse geometry. Potential mismatch of standardised Hori. (3w) is likely to happen as the estimated $\ua$ is close to its standardised value. In the bounded slack problem, slack $\slack$=0.05 is applied. We observe $\ua$ spans over a large range, suggesting that $\ua$ is not resilient to slack for the given sequence reactance. For $\radius$, its range does not overlap with $\rstd$ of Libra and Moon, so conductor type is likely to be reliably distinguished.}

\begin{table}[tbh]
\textcolor{\varcolor}{
    \setlength{\tabcolsep}{4pt}
    \caption{Inverse estimation for 3-wire geometries. All units in [mm].}
    \label{tab_inverse_OH}
    \centering
    \begin{tabular}{l l l l}
    \toprule    
    inv. geo. & inv.  & Bound-tightening & Slack bounded, $\slack$=0.05 \\
     & $\uastd$ &  $\ua$ \& $\radius$ &  Range of $\ua$ \& $\radius$ \\
    \midrule
    Hori. (3w)      & 1100 & 1155 \& 1.875 &  (729, 1500) \& (1.587, 2.017) \\
    Tri. ($\theta$=49.27°)   & 508  & 869 \& 1.875 &  (548, 1254) \& (1.587, 2.017) \\
    Tri. ($\theta$=21.67°)	& 1100 & 1100 \& 1.875 &  (694, 1500) \& (1.587, 2.017) \\
    \bottomrule
    \end{tabular}
}
\end{table}

\textcolor{\varcolor}{
Next, for cables, conductors differ in material (Al/ Cu) and number of strands. Forward calculation is  performed to understand potential mismatch, as shown in Table \ref{forward_Cable}. Again, 3- and 4-core cables are likely to be reliably distinguished due to a large difference in
zero sequence impedance. However, cables with similar ac resistance due to difference $\N$ and material type can yield similar sequence impedance.
However in reality, conductors with a large area usually come with larger number of strands $\N$. Therefore, mismatch of $\N$ may not happen often once standardised cables are considered.
}

\begin{table}[tbh]
\textcolor{\varcolor}{
    \setlength{\tabcolsep}{3pt}
    \caption{Forward calculation of make-up cables at $\T$=75\,°C, $\tnom$=1.35\,mm.}
    \label{forward_Cable}
    \centering
    \begin{tabular}{l l l l l}
    \toprule
Sensitivity  &	Cable & material & Area & $\Rzeroref,\Xzeroref,\Roneref,\Xoneref$ \\
    study            &       &          & [mm$^2$] & [$\Omega/km$]   \\
    \midrule
    base case & 3w7N & Al-1350 & 50 & (0.8395,2.2066,0.6915,0.0801) \\
    No. of strand & 3w19N &Al-1350 & 50 & (0.8395,2.202,0.6915,0.0772) \\
    Similar $\rac$ & 3w7N & Cu & 30 & (0.8645,2.2466,0.7165,0.0842) \\
    No. of core & 4w7N & Al-1350 & 50  & (1.6289,1.071,0.6916,0.0873) \\
    \bottomrule
    \end{tabular}
}
\end{table}
\textcolor{\varcolor}{
A mismatch study using standardised 4-core cable from Table \ref{tab_EQLD_construction_code} is conducted with results shown in Table \ref{tab_inverse_Cable}. The bounded slack problem with $\slack$=0.05 is solved for each combination to investigate the range of $\radius$. Mismatch happens if $\radius$ of incorrect standardised cables falls within range. We observe potential mismatch between LVABC4x25 and UGC16x4Cu, also LVABC4x95 and UGC50x4Cu\footnote{\textcolor{\varcolor}{Symbol $X$ in Table \ref{tab_inverse_Cable} may not be symmetric as it is an inverse problem with bounds stemming from domain knowledge.}}. However, they may be distinguished if utilities know whether it is below ground or not, as copper cables are likely to be employed underground.}
\begin{table}[tbh]
\textcolor{\varcolor}{
    \setlength{\tabcolsep}{3pt}
    \caption{Mismatch study of 4-core cables. Diagonal elements are correct matches ($\checkmark$). Potential mismatches are denoted by $X$.  }
    \label{tab_inverse_Cable}
    \centering
    \begin{tabular}{l l l l l l l}
    \toprule
        For. \,\, \textbackslash \,\, inv. & LVABC & LVABC & LVABC & UGC & UGC & UGC \\
    & 4x25 & 4x50 & 4x95 & 16x4Cu & 50x4Cu & 240x4Al \\
    \midrule
LVABC4x25 & $\checkmark$ & -- & -- & $X$ & -- & -- \\ 
LVABC4x50 & -- & $\checkmark$ & -- & -- & -- & -- \\ 
LVABC4x95 & -- & -- & $\checkmark$ & -- & $X$ & -- \\ 
UGC16x4Cu & $X$ & -- & -- & $\checkmark$ & -- & -- \\ 
UGC50x4Cu & -- & -- & -- & -- & $\checkmark$ & -- \\ 
UGC240x4Al & -- & -- & -- & -- & -- & $\checkmark$ \\ 
    \bottomrule
    \end{tabular}}
\end{table}

\textcolor{\varcolor}{
For the recovery of reference height $\vref$ when $\Bzeroref,\Boneref$ are provided, we look at correct combinations and investigate the recovery of $\vref$. In the feasibility problem, zero $\BZdiff$ is observed as correct combinations are used. Also $\vref$ is correctly recovered with local uniqueness observed in the bound-tightening problem. However, we find that $\vref$ varies by 2.5\,m for each \% of slack introduced, suggesting that $\vref$ is not resilient to slack in $\Bzeroref,\Boneref$.}

\subsection{Larger Test Cases}
Samples of overhead lines and cables\footnote{\textcolor{\varcolor}{In forward calculation, we assume $\tnom$ =1.5\,mm and $\vref$=-1000\,mm \cite{Ergon_underground}.}} are created to generate a variety of sequence impedance components in forward calculation, through iteration of area with step 5 mm$^2$, temperature with step 5\,°C, material in Table \ref{tab_material_parameters} and configuration in Table \ref{tab_config_parameters} over practical range as summarised in Table \ref{tab_iteration}. Shunt susceptance is excluded because it is usually not provided by utilities. 
Note that forward calculations using too small or large areas result in strand radius smaller than $\rmin$ or larger than $\rmax$ shown in Table \ref{tab_bound_value_code} are removed to prevent violation of radius bound in \eqref{eq_r_bound}.

\begin{table}[tbh]
\caption{Bound of iteration, top for OH lines and bottom for cables.}
\label{tab_iteration}
    \centering
    \begin{tabular}{l l l p{4.8cm}}
        \toprule
        Var. & Range/Set & Unit & Justification \\
        \midrule
        $\area$ & [15,240] &mm$^2$ & Min. and max. area in Table \ref{tab_EQLD_construction_code}.\\
        $\T$ & [20,75] &°C & Dc resistance temp. and layout temp. \cite{Ergon_overhead} \\
        $\material$ & $\{$Al-1350$\}$ & N/A & Assumed material for OH line. \\
        \midrule
        $\area$ & [15,240] &mm$^2$ & $\Cableconfig \in \setconfig^\text{cable}$ excepts Cable 4w48N sect.\\
        $\area$ & [185,300] &mm$^2$ & Sector area of $\N$ =37 \& $\N$ =61 \cite{sector_catalogue_NZ}.\\
        $\T$ & [20,90] &°C & Cable max. operation temp. is 90°C \cite{Ergon_underground}. \\
        $\material$ & $\setmaterial$ & N/A & Both Al and Cu are common for cables.\\
\bottomrule
    \end{tabular}
\end{table}

\subsubsection{Feasibility Problem Result}
\begin{figure}[tbh]
\centering
\includegraphics[trim={0 0.2cm 0 0.2cm},clip,width=\linewidth]{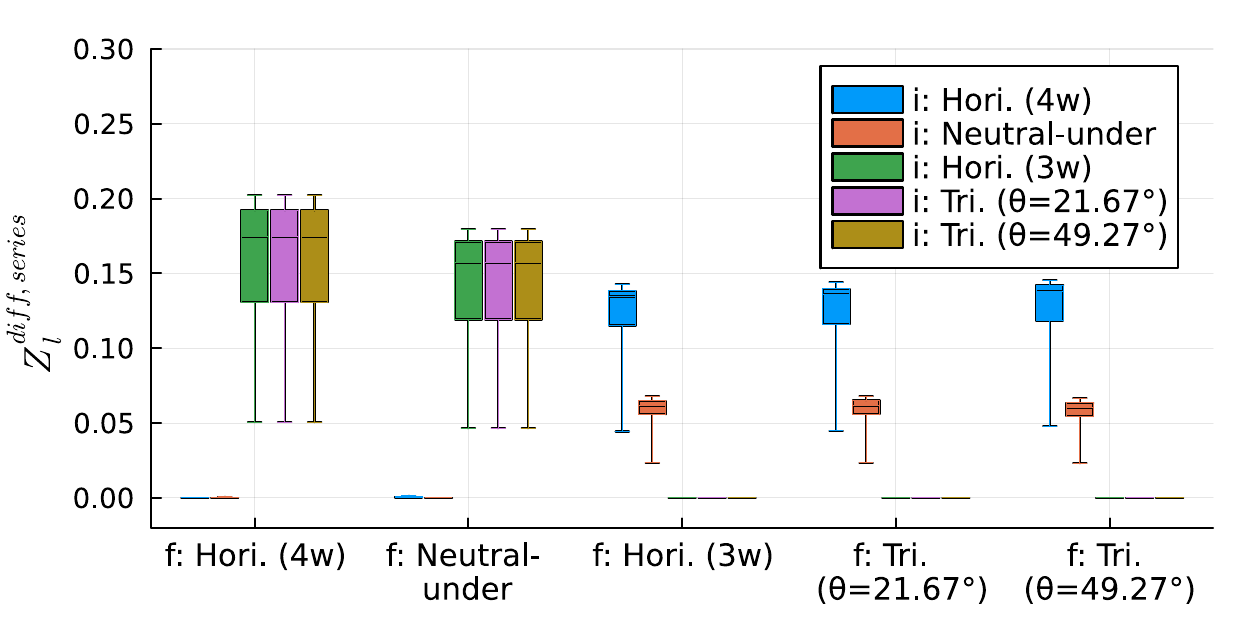}
\caption{$\Zdiff$ for overhead lines with different geometries in forward (f:) calculation  and inverse (i:) estimation.}
\label{fig:OH_Z_diff_step3}
\end{figure}

For OH lines, Fig. \ref{fig:OH_Z_diff_step3} shows that when the number of wires is incorrectly estimated in a combination, values of $\Zdiff$ range from 0.1 to 0.2. However, we observe $\Zdiff \approx$ 0 for correct $\n$. This suggests that mismatch of number of wires causes impedance highly deviated from the given sequence impedance.

\begin{figure}[tbh]
\centering
\includegraphics[trim={0 2.5cm 0 2.2cm},clip,width=\linewidth]{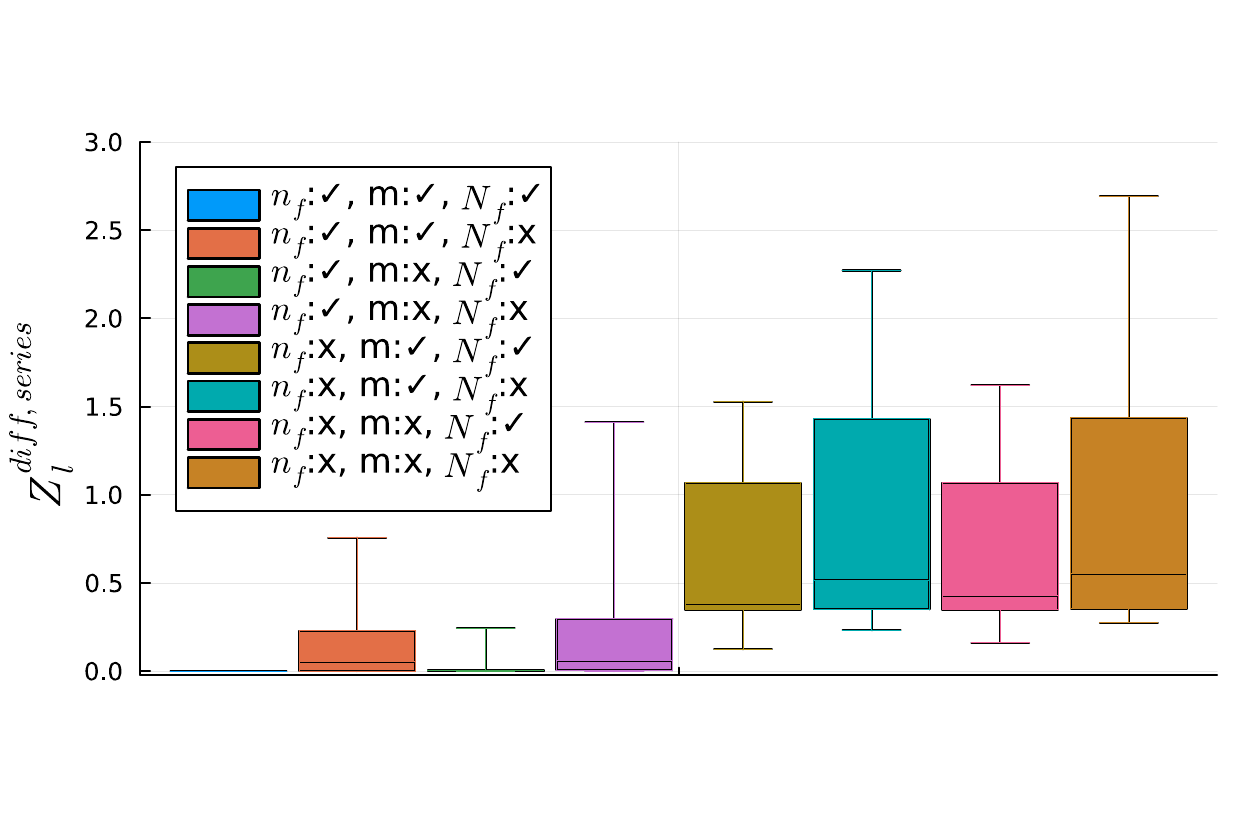}
\caption{$\Zdiff$ for cables with match and mismatch of combination of different discrete properties, as denoted by $\checkmark$ and $X$ respectively.}
\label{fig:Zdiff_truth_table}
\end{figure}

For cables, the match and mismatch of three discrete properties $\n$, $\N$ and $\material$ yield 2$^3$=8 \emph{combinations}. Fig.  \ref{fig:Zdiff_truth_table} shows that highest $\Zdiff$ occurs when mismatch of number of cores $\n$ occurs. 

\subsubsection{Bound Tightening Problem Result}
In the bound-tightening problem, (local) uniqueness over decision variables $\setvar=\{\radius,\T,\freevar{u}{\config,1} ,\freevar{u}{\config,2} ,\freevar{v}{\config,1} \}$ is studied. To prevent over-constraining, only some feasible solution $\Rzerostar, \Ronestar, \Xzerostar, \Xonestar$ are used by the bound-tightening problem, as denoted by $\checkmark$ in Table \ref{tab_DOF_bound_tightening}. 
All decision variables $\var \in \setvar$ are locally unique as validated by $\arg\min \var = \arg\max \var$ for all combinations $\combination$ up to \textcolor{\varcolor}{a tolerance of 2\,$^\circ$C for $\T$, 0.005\,mm for $\radius$, 0.04\,mm for geometry variables $\ua,\ub,\va$. For correct combinations, a maximum deviation of 0.001\% is observed for OH geometry variables over the standardised geometry parameters. Similarly for strand radius, a maximum deviation of 0.1\% is observed, which is less than the 1\% manufacturer tolerance for conductor strand radius \cite{AS3607}.}   
\begin{table}[tbh]
\setlength{\tabcolsep}{4pt}
\caption{ \textcolor{\varcolor}{Degree-of-freedom (DOF)} for different cases. Temperature and radius $\radius$ together constitutes 2 degree-of-freedom (DOF). Number of Geometry DOF in $\{ \freevar{u}{\config,1},\freevar{u}{\config,2},\freevar{v}{\config,1} \}$ depends on cases. }
\label{tab_DOF_bound_tightening}
    \centering
    \begin{tabular}{l l l l l l l l l}
        \toprule
     $|\mathcal{W}| $ &  Type   & $\T \, \& \, \radius$        & Geo. & Sum & 
      $\Rzerostar$    & $\Ronestar$    & $\Xzerostar$& $\Xonestar$ \\
                      &         & DOF    & DOF  & DOF &       &       &       &    \\
        \midrule
           3          &   OH      & 2        &  1    &  3  &       &$\tick$&$\tick$&$\tick$\\
           3          & Cable     & 2        &  1    &  3  &       &$\tick$&$\tick$&$\tick$\\
           4          &   OH      & 2        &  2    &  4  &$\tick$&$\tick$&$\tick$&$\tick$\\
           4          & Cable     & 2        &  1    &  3  &       $\tick$&$\tick$&$\tick$&$\tick$\\
        \bottomrule
    \end{tabular}
\end{table}

\subsubsection{Bounded Slack Problem Result}
\begin{figure}[tb]
\centering
\includegraphics[trim={0 0.5cm 0 0.8cm},clip,width=\linewidth]{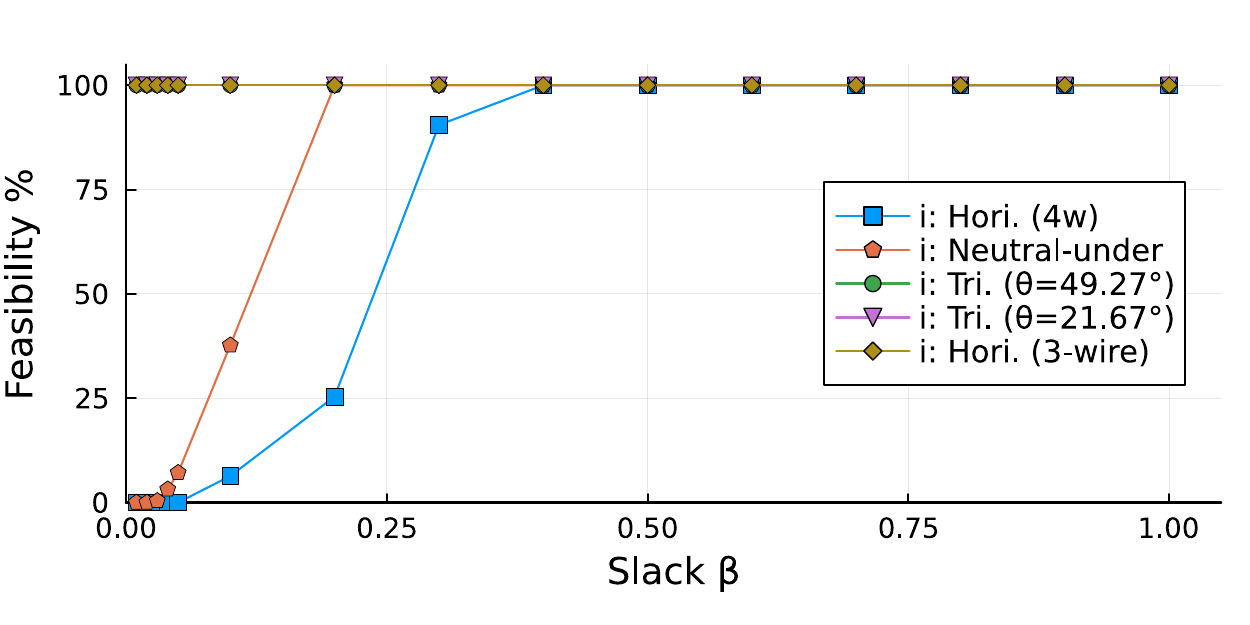}
\caption{Feasibility percentage for different inverse geometries (i:) when forward geometry is Tri. ($\theta$=21.67°).}
\label{fig:OH_massive_noise_step3_delta}
\end{figure}

For OH lines, the number of wires can be reliably distinguished.
Fig. \ref{fig:OH_massive_noise_step3_delta} shows the feasibility percentage over different inverse configuration when the forward geometry is a 3-wire Tri.  ($\theta$=21.67°). 
Inverse geometry of 4 wires are infeasible unless large amount of slack is applied, which suggests that mismatch of number of wires is not likely to happen. Nevertheless, other 3-wire inverse geometries are always feasible regardless of amount of slack because off-diagonal symmetric components which contains unbalance information are missing.

\begin{figure}[tbh]
\centering
\includegraphics[trim={0 2.2cm 0 2.5cm},clip,width=\linewidth]{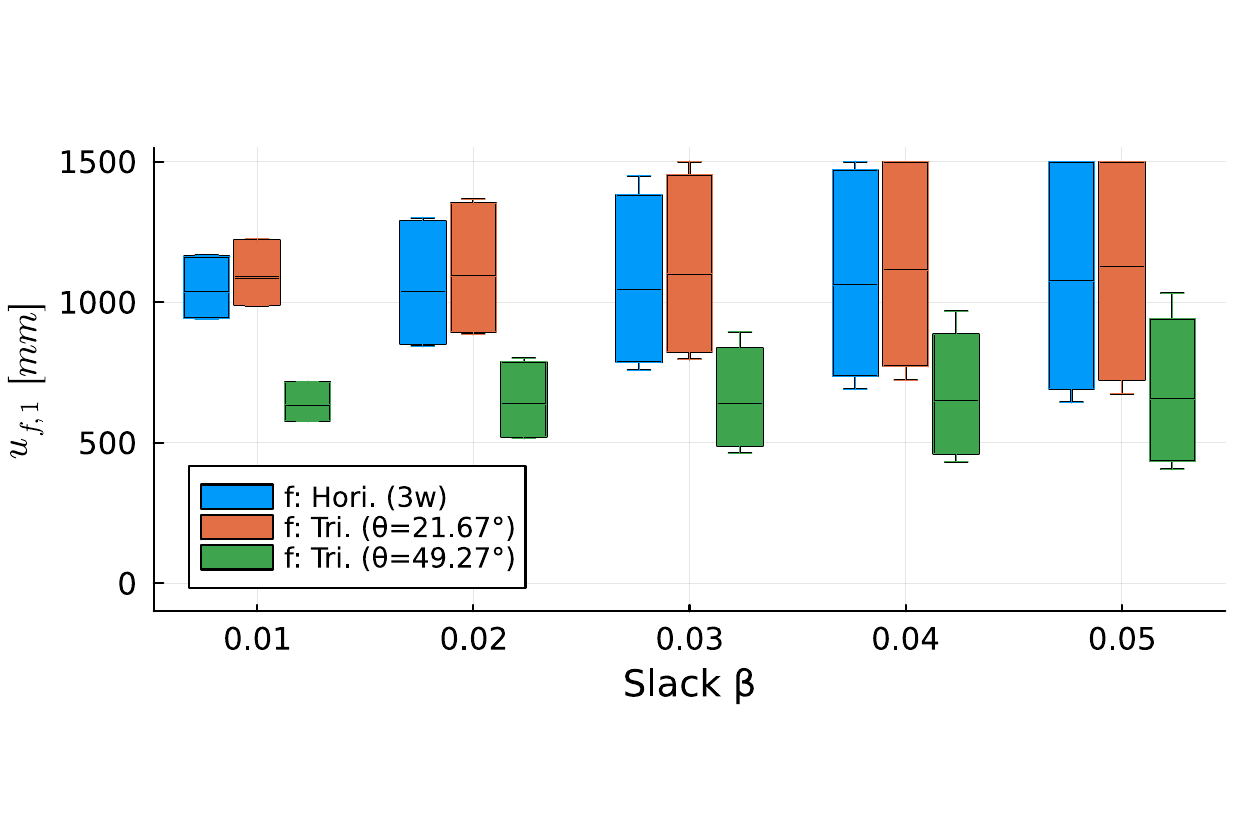}
\caption{Range of $\freevar{u}{\config,1}$ with inverse geometry Tri. ($\theta$=21.67$^\circ$) over 3-wire forward (f:) geometries. Variable $\freevar{u}{\config,1}$ has an upper bound of $\umaxOH$=1500\,mm.}
\label{fig:u1_delta_step3}
\end{figure}

\textcolor{\varcolor}{The ability to distinguish geometry types is now investigated. Fig. \ref{fig:u1_delta_step3} shows the range of $\ua$ under different slacks for all 3-wire forward geometries in the bounded slack problem. All boxes in the box-plot have a large range of $\ua$ even under small slack, which suggests that geometry variables are highly sensitive to slack in given sequence values. Next, overlapping of  $\ua$ happens between Hori. (3w) and Tri. ($\theta$=21.67°) at all slack. However, Tri. ($\theta$=49.27°) can be reliably distinguished or eliminated assuming up to 3\% in given sequence components. This aligns with the insight from previous study that similar standardized geometries with similar sequence reactance are hard to distinguish. 
}   

\begin{figure}[tbh]
\centering
\includegraphics[trim={0 2cm 0 2cm},clip,width=\linewidth]{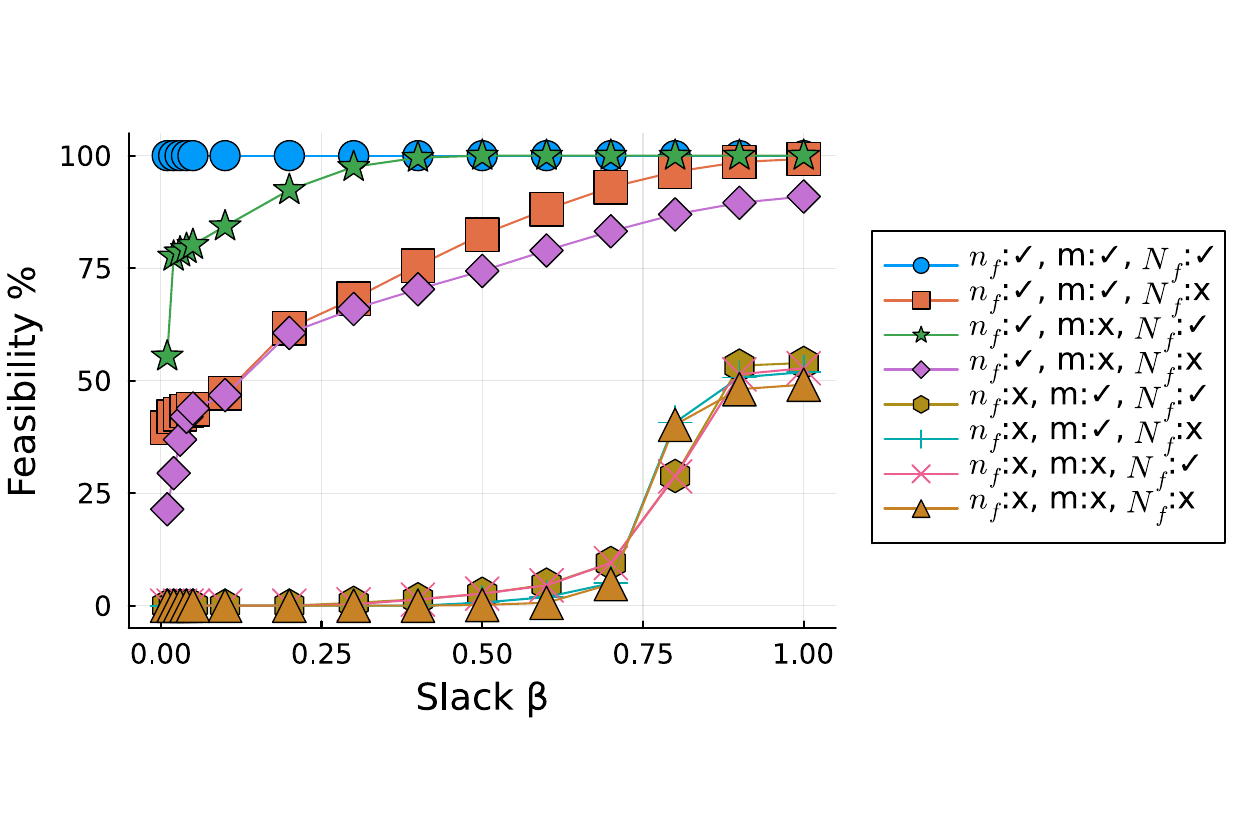}
\caption{Match and mismatch of combination of discrete properties against forward calculation are denoted by $\checkmark$ and $X$ respectively.}
\label{fig:percentage_truth_table}
\end{figure}

For cable, Fig. \ref{fig:percentage_truth_table} shows the overall feasibility percentage for different \emph{combinations} of match or mismatch of discrete properties. It shows mismatch of number of core makes a combination to become infeasible unless high slack is applied. Mismatch of number of strands is more likely to cause an combination to be infeasible than mismatch of material. 

\subsection{\textcolor{\varcolor}{Real World Examples}}
\textcolor{\varcolor}{
We now demonstrate the recovery using example impedance data \cite{Ergon_overhead, Ergon_underground} from an Australian utility and the result is summarised in Table \ref{tab_EQLD_forward_comparison}. Both diagonal sequence components
and line information are provided. However, the exact line parameters, e.g. $\uastd$, are not provided. For overhead lines, the forward calculation result matches well with the provided sequence component values. For cables, it is clearly seen that zero sequence impedance is derived from positive ones, with $\Rzeroref$=$4\cdot\Roneref$ and $\Xzeroref$=$\Xoneref$. This highlights the need for data validation specific to sequence impedance values.
}
\begin{table}[tbh]
\textcolor{\varcolor}{
    \setlength{\tabcolsep}{2pt}
    \caption{Cable at 20°C, Overhead lines Mars at 75°C}
    \label{tab_EQLD_forward_comparison}
    \centering
    \begin{tabular}{l l l}
    \toprule
    Line & ($\Rzeroref,\Xzeroref,\Roneref,\Xoneref$) & ($\Rzeroref,\Xzeroref,\Roneref,\Xoneref$)\\
    code & utility [$\Omega/km$] & forward cal. [$\Omega/km$] \\    
    \midrule
    Tri. (21.67°) & (0.600,1.613,0.452,0.356) & (0.5952,1.5873,0.4472,0.3692) \\
    Tri. (49.27°) & (0.600,1.631,0.452,0.347) & (0.5952,1.6547,0.4472,0.3355) \\
    UGC16x4Cu & (4.6,0.089,1.15,0.089) & (2.096,1.5195,1.1185,0.0917) \\
    UGC50x4Cu & (1.55,0.082,0.388,0.082) & (1.0792,0.6336,0.369,0.0891) \\
    UGC240x4Al & (0.500,0.062,0.126,0.062) & (0.4021,0.2933,0.1183,0.0676) \\
    \bottomrule
    \end{tabular}
    }
\end{table}

\textcolor{\varcolor}{
The percentage mismatch over $\rstd$ for cables and $\uastd$ for overhead lines are summarised in Table \ref{tab_EQLD_recovery}. Lowest percentage is chosen as the candidate for inverse estimation. For overhead lines, only 3-wire combinations have $\Zdiff$=0 so 4-wire combinations are eliminated.
Next, Tri. (49.27°) is matched correctly but 20\% difference over $\uastd$ is still shown. It is expected because sequence reactance is not sensitive to $\ua$. For Mars, Tri. (21.67°), mismatch happens as the best candidate is Hori. (3w). However, mismatch risk is low because they share close sequence impedance. Another candidate Tri. (49.27°) with higher mismatch risk (more distinct reactance) is eliminated due to higher \% mismatch for $\ua$. In both cases, recovered $\radius$ is 1.88\,mm and hence conductor Mars ($\rstd$=1.875\,mm) is successfully recovered.}

\textcolor{\varcolor}{
For cables, if the provided zero sequence values are used, $\Zdiff \geq$ 0.25 is observed for all combinations, meaning an average of 25\% deviation for all sequence components from the Carson-based impedance model. This suggests that no combination can well explain the data, and data validation is needed. Therefore, the provided zero sequence values are deemed missing now for the recovery process. Without the provided zero sequence values, 3-core combinations can have zero $\Zdiff$ even for a 4-core cable. However, as Australian utilities generally know the number of phases and neutral in LV networks, 3-core combinations can be eliminated. Next, for the remaining 4-core combinations, if $\Zdiff \neq$ 0, the associated candidates are eliminated. Similar to the analysis for Table \ref{tab_inverse_Cable}, potential mismatch happens for UGC16x4Cu and UGC50x4Cu in Table \ref{tab_EQLD_recovery} with low \% mismatch of $\rstd$ for wrong candidates. If temperature information is given, UGC50x4Cu can be reliably distinguished with low \% mismatch compared with others. For UGC16x4Cu, we need to know its material type or whether it is buried underground for reliable distinction.}
\begin{table}[tbh]
\textcolor{\varcolor}{
    \setlength{\tabcolsep}{2.5pt}
    \caption{Percentage difference of $\ua$ for standardized OH line and $\radius$ for standardized cable. ``[]'' indicates when $\T$ is known whereas ``\textbackslash'' candidates with $\Zdiff \neq$ 0 for its combination. }
    \label{tab_EQLD_recovery}
    \centering
    \begin{tabular}{p{2.7 cm} p{1.8cm} p{1.8cm} p{1.8cm} }
    \toprule
    utility data\textbackslash inv. geo & {Hori. (3w)} & {Tri. (49.27°)} & {Tri. (21.67°)} \\
    \midrule
    Mars, Tri. (49.27°) & {26.0} & {20.5} & {29.5} \\
    Mars, Tri. (21.67°) & {14.6} & {39.2 } & {18.6 } \\
    \bottomrule
    \end{tabular}
    \begin{tabular}{l l l l l l l}
    utility \textbackslash inv. & LVABC & LVABC & LVABC & UGC & UGC & UGC \\
    data & 4x25 & 4x50 & 4x95 & 16x4Cu & 50x4Cu & 240x4Al \\
    \midrule
UGC16x4Cu & 7.1[3.8] & 20.4 & \textbackslash & 13.5[0.0] & 34.8 & \textbackslash \\ 
UGC50x4Cu & 80.0 & 33.8 & 4.9[12.3] & 95.4 & 12.2[2.5] & \textbackslash \\ 
UGC240x4Al & \textbackslash & \textbackslash & \textbackslash & \textbackslash & \textbackslash & 3.6 \\
    \bottomrule
    \end{tabular}
    }
\end{table}

\section{Conclusions} \label{sec_conclusions}
\textcolor{\varcolor}{
This paper formulates the recovery of low-voltage network data as an inverse problem from first principles - i.e. Carson's equations -  determining the best fit w.r.t. the physics for untranposed overhead lines and cables. 
The inverse problem methodology involves three optimisation-based problems. 
In the feasibility problem, the number of conductors can be reliably distinguished. In the bound-tightening problem, we observed local uniqueness for the recovered conductor and geometry variables, which can be correctly recovered for matched material and configuration types. In the bounded slack problem, for overhead lines, distinct geometry types can be distinguished assuming that limited uncertainty presents in sequence reactance. Next, standardised cables with distinct ac resistance can be distinguished. Mismatch happens for overhead lines with similar geometry types and cables with similar ac resistance. However, mismatch risk is low because they share similar impedance. We also demonstrate the recovery of missing data through a case study using utility data, which highlights the importance of data validation. Exploiting the structure of Carson's equations when learning line impedance matrices from smart meter data is a direction of future work. 
}


\bibliographystyle{IEEEtran}
\bibliography{library} 
\end{document}